\documentclass[11pt,a4paper,reqno]{amsart}

\newcommand{\ol}[1]{\overline{#1}}

\usepackage{tikz}
\usetikzlibrary{decorations.pathreplacing,shapes,arrows}
\newcommand{\arrowdraw}[2]{\draw [postaction={decorate}] (#1)--(#2);}
\newcommand{\darrowdraw}[2]{\draw [dashed, postaction={decorate}] (#1)--(#2);}

\usetikzlibrary{decorations.markings}
\usetikzlibrary{arrows,matrix}
\usetikzlibrary{snakes}
\usepgflibrary{arrows}
\tikzset{
    >=stealth',
    punkt/.style={
           rectangle,
           rounded corners,
           draw=black, very thick,
           text width=6.5em,
           minimum height=2em,
           text centered},
    pil/.style={
           ->,
           thick,
           shorten <=2pt,
           shorten >=2pt,}
}
\tikzset{every loop/.style={min distance=2mm,in=225,out=135,looseness=10}}

\usepackage{tikz}
\usetikzlibrary{decorations.markings}
\usetikzlibrary{arrows,matrix}
\usepgflibrary{arrows}

\newcommand\GreenL{\mathcal{L}}
\newcommand\GreenR{\mathcal{R}}

\newcommand\GreenJ{\mathcal{J}}

\usepackage{clrscode}
\usepackage{mathrsfs,amssymb}

\usepackage{amsmath, amsthm}

\newtheorem{theorem}{Theorem}[section]

\newtheorem{question}[theorem]{Question}
\newtheorem{lemma}[theorem]{Lemma}
\newtheorem{corollary}[theorem]{Corollary}
\newtheorem{conjecture}[theorem]{Conjecture}
\newtheorem{proposition}[theorem]{Proposition}

\theoremstyle{definition}
\newtheorem{example}[theorem]{Example}
\newtheorem{remark}[theorem]{Remark}
\newtheorem{definition}[theorem]{Definition}

\begin{document}

\title[Cogrowth, Amenability and Spectral Radius]{On Cogrowth, Amenability and the Spectral Radius of a Random Walk on a Semigroup}

\maketitle

\begin{center}
    ROBERT D. GRAY\footnote{School of Mathematics, University of East Anglia, Norwich NR4 7TJ, England.
Email \texttt{Robert.D.Gray@uea.ac.uk}.}
\ and \ MARK KAMBITES\footnote{School of Mathematics, University of Manchester, Manchester M13 9PL, England. Email \texttt{Mark.Kambites@manchester.ac.uk}.

This research was supported by the EPSRC-funded project EP/N033353/1 `Special inverse monoids: subgroups, structure, geometry, rewriting systems and the word problem'.
} \\
\end{center}

\begin{abstract}
We introduce two natural notions of cogrowth for finitely generated semigroups --- one \textit{local} and one \textit{global} --- and
study their relationship with amenability and random walks. We establish the minimal and maximal possible values for cogrowth
rates, and show that non-monogenic free semigroups are exactly characterised by minimal global cogrowth. We consider the relationship with
cogrowth for groups and with amenability of semigroups.  We also study the relationship with random walks on finitely generated semigroups, and in particular the spectral radius of the associated Markov operators (when defined) on $\ell_2$-spaces. We show that either of maximal global cogrowth or the weak F\o lner condition suffices for its spectral radius to be at least $1$; since left amenability implies the weak F\o lner condition, this represents a generalisation to semigroups of one implication of Kesten's Theorem for groups. 
By combining with known results about amenability, we are able to establish a number of new sufficient conditions for (left or right) amenability in broad classes of semigroups. In particular, maximal local cogrowth left implies amenability in any left reversible semigroup, while maximal global cogrowth (which is a much weaker property) suffices for left amenability in an extremely broad class of semigroups encompassing all inverse semigroups, left reversible left cancellative semigroups and left reversible regular semigroups.
\end{abstract}

\section{Introduction}\label{sec_intro}

The concept of \textit{cogrowth} forms a natural bridge between combinatorial, probabilistic and analytic approaches to the study of finitely
generated groups. Defined in terms of generators, as an elementary asymptotic invariant in a natural sense dual to growth, the cogrowth
of a group has a natural interpretation in terms of divergence rates of a random walk on the Cayley graph, and can be used to characterise amenability \cite{Cohen82,Grigorchuk80} of the group.

The areas of group theory touched upon by the study of cogrowth all have counterparts in semigroup theory. There has long been interest, both intrinsic and motivated by applications in computer science, in the combinatorial study of (usually finitely generated) semigroups using generators and relations (see for example \cite{Birget1998, Guba1997, Kobayashi2005, Robertson1998, Steinberg2003}).
Random walks on semigroups are of longstanding interest, and have recently received additional attention due to connections with representation theory \cite{Bidigare1999, Brown2000, Brown1998, Margolis2015}). The focus so far has been chiefly on the finite case (most notably on face monoids of hyperplane arrangements), but the extension to finitely generated semigroups is arguably the next natural step. 
Finally, amenability for semigroups has been an active area of research from the 1950s \cite{Day57} to the present day; see 
\cite{Barber15,Bartoszek17,Bergelson12,Cecc15,Hindman09} for some examples of recent developments and applications in
other areas and \cite{GrayKambitesAmenable,Takahashi03} for advances relating specifically to finitely generated semigroups.  It therefore seems apposite to ask whether whether the concept of cogrowth makes sense in a semigroup setting, and if so whether it is capable of forming a similar bridge between these areas.

In this paper, we introduce, and begin the study of, two natural notions of cogrowth for finitely generated semigroups and monoids --- one \textit{local} and one \textit{global} ---
and study their relationship with random walks and with amenability.
\textit{Local} cogrowth is the most immediately obvious extension to semigroups of the usual definition of cogrowth in groups; we shall see that for semigroups in general it is not capable of capturing the same level of structural and dynamic information as in groups, but it is interesting in some cases and also important as a tool for studying global cogrowth. \textit{Global} cogrowth is a slightly more subtle definition which for semigroups turns out to contain much more information, while for groups remaining essentially equivalent to the usual definition; we contend that this is the ``correct'' notion of cogrowth for semigroups. 

We establish the minimal and maximal possible values for cogrowth
rates, and show that (excepting the $1$-generator case, which has a pathological aspect) free semigroups are exactly characterised by minimal global cogrowth (Theorem~\ref{thm_free}). Where the semigroup happens to be a group with a symmetric generating set, our notions are closely related to each other and to the group-theoretic notion of cogrowth; in particular, we deduce from the Grigorchuk-Cohen cogrowth theorem \cite{Cohen82,Grigorchuk80} that maximality (of either) characterises amenability (Theorem~\ref{thm_groupamen}). In greater generality (for semigroups, or even for groups with non-symmetric generating sets) we show that maximality of local cogrowth is not invariant under change of generators, so cannot be expected to characterise any abstract property of semigroups (Example~\ref{example_groupnonsymmetric_local}). It remains unclear whether
maximality of global cogrowth is invariant under change of generators, but we are able to show that if a semigroup 
has a generating set with respect to which it has maximal global cogrowth, then every finite
subset is contained in such a generating set (Theorem~\ref{thm_supersets}).

To establish the connection with random walks, we follow Kesten's approach \cite{Kesten59b,Kesten59a} of studying the spectral radius of the associated Markov operator on the $\ell_2$-space of the semigroup. Unlike in groups, this operator is not always defined; Proposition~\ref{prop_operatorcondition} below gives an exact algebraic description of when it is defined. When it is defined, it is bounded and we show (Theorems~\ref{thm_kestenspectralradius} and \ref{thm_fcimpspecrad}) that any of maximal local cogrowth, maximal global cogrowth, the weak F\o lner condition or left amenability suffices for its
spectral radius to be at least $1$. The latter criterion generalises to semigroups one implication of Kesten's
Theorem \cite{Kesten59b} for groups. 

By combining our results with known results about amenability, we are able to deduce a number of new sufficient conditions for left amenability in left reversible semigroups (left reversibility being an elementary structural precondition for left amenability).  Maximal local cogrowth (at any element) suffices for left amenability in any left reversible semigroup (Theorem~\ref{thm_localimpamenable}). The weaker condition of maximal global cogrowth suffices for left amenability in left reversible semigroups satifying the \textit{Klawe condition} (Theorem~\ref{thm_cogrowthamenable_Klawe}), and hence in the extreme broad class of \textit{near left cancellative} semigroups introduced in \cite{GrayKambitesAmenable}. Maximal global cogrowth also suffices for left amenability 
in left reversible semigroups whose maximal right cancellative quotient has a minimal ideal (Theorem~\ref{thm_cogrowthamenable}). Both of
these results apply in particular to inverse semigroups, allowing us to show that a finitely generated inverse semigroup of maximal global cogrowth
is necessarily amenable (which can be deduced from known results when the generating set is symmetric, but is a new result where it is not). The main known and unknown implications between the various conditions considered are summarised diagramatically in Figure~\ref{fig_relationships} at the end of the paper.

\section{Preliminaries}

Throughout this paper we assume basic familiarity with some elementary concepts of semigroup theory, such as idempotents, ideals, Green's equivalence relations and pre-orders, (von Neumann) regular elements and inverse semigroups; for an introduction to these the reader is directed to any of the standard texts, such as \cite{Grillet95,Howie95}. For the avoidance of confusion, we recall that a \textit{simple} semigroup or monoid is one without proper 2-sided ideals; this is not the same thing as having no proper homomorphic images, semigroups and monoids with the latter property being termed \textit{congruence-free}. A \textit{subgroup} of a semigroup means a subsemigroup which forms a group under the inherited multiplication (the identity element of which may be any idempotent of the semigroup); a \textit{maximal subgroup} is a subgroup not properly contained in another subgroup; there is no requirement that it be a proper subsemigroup, so if the semigroup happens to be a group then, in contrast with the usual terminology in group theory, it is its own unique maximal subgroup. If $S$ is a semigroup without an identity we write $S^1$ to denote the
monoid obtained by adjoining an identity element; if $S$ is already a monoid then $S^1$ is defined to be $S$.

While the earlier sections should be intelligible to semigroup theorists with no background in analysis or amenability, Section~\ref{sec_randomwalk} is concerned with random walks and associated Markov operators, while Section~\ref{sec_amenable} involves amenability of semigroups. Very brief introductions to each are given in the respective sections, but for a more leisurely and comprehensive introduction we direct the reader to expositions such as \cite{Cecc10, Paterson88} for amenability of semigroups and \cite{Woess00} for concepts related to random walks.
 
In the rest of this section we fix some notation and terminology. If $X$ is a set of symbols, we write $X^*$ (respectively, $X^+$) to represent respectively the free monoid (free semigroup)
on $X$, that is, the set of finite words (respectively, finite non-empty words) over $X$ under the operation of concatenation.
For $u \in X^*$ we write $|u|$ for the length of $u$.

A (\textit{finite}) \textit{choice of generators} for a semigroup (or monoid) $S$ is
a (finite) set $X$ of symbols together with a surjective morphism from the free
semigroup $X^+$ (or the free monoid $X^*$) onto $S$. For
words $u$ and $v$ in $X^+$ or $X^*$ we
write $\ol{u}$ for the image of $u$ in $S$, and
$u \equiv v$ to mean that $\ol{u} = \ol{v}$ in $S$.
The canonical example is where $X$ is a generating subset of
$S$ with the map given by multiplication in the obvious way,
but our definition is slightly more general because multiple
elements of $X$ may map to the same element of $S$,
allowing us to consider generators ``with multiplicity''.
We shall often simply say that $X$ is a choice of generators for $S$,
leaving the map implicit. By a \textit{finitely generated semigroup} (\textit{monoid})
we mean a semigroup (monoid) equipped with a fixed finite choice of generators.
When $S$ is a group or inverse semigroup, the choice of generators is called \textit{symmetric}
if for every $s \in S$ the same number (possibly $0$) of generators in $X$ represent $s$ and $s^{-1}$.
Throughout the paper, we shall assume unless stated otherwise that $S$ is a semigroup with a fixed finite choice of generators $X$.

We take $\mathbb{N}$ to be the set of positive natural numbers, excluding $0$. If $A$ is a subset of a set $B$ we write $\chi_A : B \to \mathbb{R}$ for the characteristic function of $A$; if $a \in B$ we write $\chi_a$ for $\chi_{\lbrace a \rbrace}$.

\section{Local Cogrowth}

In this section we introduce the first of two notions of cogrowth for finitely generated semigroups, which is the \textit{local cogrowth} at a given element. Although this turns out to be (in our opinion) less important than the notion of \textit{global cogrowth} to be introduced in the next section, it is closely related to global cogrowth and provides a useful tool for studying it.

\begin{definition}
The \textit{local cogrowth function} at an element $s \in S$ is the function
$\lambda^{S,X}_s : \mathbb{N} \to \mathbb{N}$ which maps $n$ to the number
of distinct words of length $n$ in $X^+$ representing $s$.

The \textit{local cogrowth rate} $\lambda^{S,X}_s$ at $s \in S$ is
$$\limsup_{n \to \infty} \ (\lambda^{S,X}_s(n))^{1/n}.$$
Where the semigroup and choice of generators are clear, we drop
the superscripts and write simply $\lambda_s(n)$ and $\lambda_s$.
\end{definition}

\begin{remark}
Clearly since there are only $|X|^n$ words of length $n$ we have
$$0 \ \leq \ \lambda_s(n) \ \leq \ |X|^n$$
for all $n$. Thus, the sequence $\lambda_s(n)^{1/n}$ is bounded and hence $\lambda_s$ is defined
and $0 \leq \lambda_s \leq |X|$.
\end{remark}

\begin{definition}
If $\lambda_s = |X|$ we say that $s$ has \textit{maximal local cogrowth}, or that $S$ has \textit{maximal local cogrowth} (\textit{at the element $s$}).
\end{definition}

\begin{remark}
Semigroups with maximal local cogrowth are examples of the \textit{$A$-semigroups} introduced and studied by Gerl \cite{Gerl73}. See Remark~\ref{remark_gerl} below for further discussion of the relationship between these two concepts.
\end{remark}

\begin{remark}
In the special case that $S$ is a group with identity $1$ and $X$ is a symmetric generating set,
$\lambda_1$ corresponds with the formulation of cogrowth considered in for example \cite{Elder12}. It differs technically from the standard formulation for groups, in that it counts all words representing the identity, rather than just freely reduced words. However, it is closely related and still suffices
to characterise amenability (see \cite{Elder12} and Theorem~\ref{thm_groupamen} below). 
\end{remark}

\begin{remark}
The local cogrowth of a semigroup is not left or right ``sided'': more precisely the local cogrowth of an element $s \in S$
is the same as its cogrowth in the \textit{opposite semigroup} $S^{op}$ with the same set of elements, order of
multiplication reversed, and obvious choice of generators. This immediately implies that one cannot hope to characterise certain properties of abstract semigroups,
such as left amenability or right amenability (see Section~\ref{sec_amenable} below), which differ in $S$ and $S^{op}$, in terms of local cogrowth.
\end{remark}

\begin{proposition}
For any $s \in S$ we have $\lambda_s = 0$ if and only if $s$ is represented by only finitely many words over $X$, and $\lambda_s \geq 1$ otherwise. 
\end{proposition}
\begin{proof}
If $s$ is represented by only finitely many words then the sequence $\lambda_s(n)$ is eventually $0$, so $\lambda_s(n)^{1/n}$ is eventually $0$.
Otherwise $\lambda_s(n) \geq 1$ for infinitely many terms, so $\lambda_s(n)^{1/n} \geq 1$ for infinitely many terms. 
\end{proof}

\begin{remark}\label{remark_localrandomwalk}
The local cogrowth at an element $s$ may be thought of as measuring the probability that a random
walk in the semigroup ends at $s$. More precisely, $\lambda_s(n) / |X|^n$ is the probability
that a random walk of length $n$, starting at the identity element of $S$ (or of $S^1$ if necessary) and multiplying
consecutively on the right (say) by $n$ generators chosen uniformly
at random from $X$, ends at $s$. Notice that $\lambda_s < |X|$ if and only if this probability decays exponentially as $n \to \infty$.
\end{remark}

\begin{proposition}\label{prop_jorder}
If $s, t \in S$ are such that $s \leq_\GreenJ t$ then $\lambda_s \geq \lambda_t$. In particular, local cogrowth is a $\GreenJ$-class invariant and if $S$
is simple then local cogrowth is constant across the semigroup.
\end{proposition}
\begin{proof}
Suppose $s \leq_\GreenJ t$. Fix words $u, v \in X^*$ such
that $\ol{u} t \ol{v} = s$, and let $m = |uv|$.
Then for each word $w \in X^n$ representing $t$ we obtain a different word $uwv \in X^{n+m}$ representing
$s$,  so we have
$$\lambda_t(n) \leq \lambda_s(n+m) \textrm{ for all $n$},$$
from which the claim easily follows.
\end{proof}
\begin{definition}
If $J$ is a $\GreenJ$-class of $S$ then we use the term \textit{local cogrowth at $J$}, and write $\lambda_J$, for the local cogrowth of any element $s \in J$.
\end{definition}

\begin{proposition}\label{prop_localminimalideal}
If $S$ has an element of maximal local cogrowth then the set of all such
elements forms a minimal ideal of $S$.
\end{proposition}
\begin{proof}
It is immediate from Proposition~\ref{prop_jorder} that the set of elements of maximal
local cogrowth (presuming it to be non-empty) forms an ideal; call it $I$. Suppose
for a contradiction that it is not minimal. Let $J$ be an ideal of $S$
properly contained in $I$,
and choose elements $j \in J$ and $i \in I \setminus J$. Let $w \in X^+$ be a word
representing $j$. Then all words containing $w$ as a factor represent elements of $J$,
so in particular no such word can represent $i$.

Hence, every word representing $i$ can be written in the form $u_1 \dots u_p v$
where $|v| < |w|$ and for $1 \leq q \leq p$ we have $|u_q| = |w|$ and $u_q \neq w$. The number of ways to choose  such a
word of some length $n$ is bounded above by
$$(|X|^{|w|}-1)^p \ |X|^{|w|} \ \leq \ |X|^{|w|} \ \left( (|X|^{|w|}-1)^{1/|w|} \right)^n$$
where $p$ is the integer part of $\frac{n}{|w|}$ and the inequality holds because $|X|^{|w|} -1 \geq 1$ (unless $|X| = 1$ in which case
both sides are $0$). Since $|X|^{|w|}$ is a constant independent of $n$, we obtain
$$\lambda_i \ = \ \limsup_{n \to \infty} \lambda_i(n)^{1/n} \ \leq \ (|X|^{|w|}-1)^{1/|w|} \ < \ |X|$$
contradicting the assumption that $S$ has maximal local cogrowth at $i$.
\end{proof}

\begin{corollary}\label{cor_maxlocalsimple}
If a monoid has maximal local cogrowth at the identity then it must
be simple (that is, have no proper $2$-sided ideals).
\end{corollary}

\begin{proposition}\label{prop_localquotient}
Let $f : S \to T$ be a surjective morphism and consider $X$ as a choice of generators
for $T$ (with the obvious map obtained by composing the map from $X^+$ to $S$ with the
map $f$ from $S$ to $T$). Then for any $s \in S$ we have
$$\lambda^{T,X}_{f(s)}(n) \ \geq \ \lambda^{S,X}_s(n), \textrm{ for all } n$$
and hence $\lambda^{T,X}_{f(s)} \geq \lambda^{S,X}_s$.
\end{proposition}
\begin{proof}
This is immediate from that fact that any word representing $s$ in 
$S$ represents $f(s)$ in $T$.
\end{proof}

\begin{example}\label{example_bicycliclocal}
Recall that the \textit{bicyclic monoid} is given by the
monoid presentation
$$\langle b, c \mid bc = 1 \rangle.$$
The bicyclic monoid is ubiquitous in infinite semigroup theory and also (as a natural algebraic
model of a counter or a one-sided shift) in many other areas of mathematics and theoretical
computer science.
Letting $X = \lbrace b, c \rbrace$ be the obvious choice of generators, there are no words of odd length representing
the identity in this monoid, while the number of words of length $2k$ representing
the identity is well-known (and easily seen) to be the Catalan number $C_k$. Thus, we have
$$\lambda_1(n) \ = \ \begin{cases} 
     0 &\textrm{ if $n$ is odd} \\
     C_{n/2} &\textrm{ if $n$ is even.}
\end{cases}
$$
As a consequence of Stirling's approximation,
$$C_k \ \approx \ \frac{4^k}{k^{3/2} \sqrt{\pi}}$$
so putting $n = 2k$ we get $\lambda_1 = 2$. Thus 
by Corollary~\ref{cor_maxlocalsimple} we obtain the (well-known) fact that the bicyclic
monoid is simple, and we can deduce by Proposition~\ref{prop_jorder} that
$\lambda_s = 2$ for all elements $s$. In other words, all
elements have maximal local cogrowth.
\end{example}

\begin{example}\label{example_freecommutativelocal}
Let $S$ be the free commutative monoid $\langle a, b \mid ab = ba \rangle$
with the obvious choice of generators.
For any element $s$ there are only finitely many words representing $s$, so we have $\lambda_s(n) = 0$ for sufficiently
large $n$ and hence $\lambda_s = 0$ for all $s$.
\end{example}

The argument in Example~\ref{example_freecommutativelocal} clearly applies also to a broad class of
finitely generated semigroups with the property that each element has
only finitely many representatives. This includes all finitely generated semigroups admitting
homogeneous presentations (that is, finite or infinite presentations in which the left-hand-side
and right-hand-side of each relation have the same length), and all semigroups
given by finite presentations satisfying the small overlap condition
$C(4)$ (see \cite{Higgins92,K_smallover1,Remmers80}) which in a precise statistical sense is
``almost all'' finitely presented semigroups \cite{K_generic}.

Although it will turn out that the local cogrowth rate of elements is interesting
in certain types of semigroups, these examples show that there are many
other interesting semigroups where it can tell us little or nothing. This is one
of the motivations for a related ``global'' definition of
cogrowth, which will be introduced in the next section.

\begin{proposition}\label{prop_minimalideal}
If $S$ has a finite minimal left or right ideal $I$, then the elements of $I$
have maximal local cogrowth. In particular, any left, right or two-sided zero
element of a semigroup has maximal local cogrowth.
\end{proposition}
\begin{proof}
We treat the case that $I$ is a finite right ideal, the case of a left ideal
being dual.
Suppose $I$ has $p$ elements. Fix a word $w \in X^+$ representing some
element of $I$. For any $n > |w|$ each
word $u \in X^{n - |w|}$ yields a word $wu \in X^n$ representing an
element of $I$. Thus, for each $n$, there are at least $|X|^{n - |w|}$
such words,
which by the pigeon hole principle means that some element of $I$
has at least
$$\frac{1}{p} |X|^{n - |w|} \ = \ \frac{1}{p |X|^{|w|}} |X|^n$$
representatives of length $n$ for infinitely many $n$. Clearly,
this element will have local cogrowth $|X|$. Since $I$ is a minimal
right ideal it is contained in a  $\GreenJ$-class so it follows by
Proposition~\ref{prop_jorder} that every element
of $I$ has local cogrowth $|X|$.
\end{proof}
\begin{corollary}\label{cor_finitelocal}
Every finite semigroup has maximal local cogrowth.
\end{corollary}

Propositions~\ref{prop_localminimalideal} and \ref{prop_minimalideal} give respectively a necessary condition and a sufficient condition
for the existence of elements of maximal local cogrowth, in terms of the ideal structure of the semigroup. One might hope that they are
converging upon an abstract, purely algebraic characterisation of this property, without reference to generators. However, the following example
shows that this property is actually dependent upon the choice of generators, so such a characterisation cannot exist:

\begin{example}\label{example_groupnonsymmetric_local}
Let $G = \mathbb{Z}$ with choice of generators $X = \lbrace a,b \rbrace$ where $a$ represents $+1$ and $b$ represents $-1$. We
shall compute the local cogrowth function at the identity element: $\lambda_0(n)$.
Clearly only words of even length can represent the identity so
$\lambda^{G,X}_0(n) = 0$ for $n$ odd. 
A word of even length $n = 2k$ represents $0$ if and only if it contains $b$ in exactly $k$ positions. There are $\binom{2k}{k}$ ways to choose the positions of the $b$s and the remaining positions must contain $a$s,
so using Stirling's approximation:
$$\lambda^{G,X}_0(n) \ = \ \binom{2k}{k}  \ \approx \ \frac{4^k}{\sqrt{\pi k}} \ = \ \frac{1}{\sqrt{\pi \frac{n}{2}}} 2^n$$
whence $\lambda_0^{G,X} = 2$, which is maximal. Since $\mathbb{Z}$ (like all groups) is a simple semigroup, it follows by Proposition~\ref{prop_jorder} that $\lambda_z = 2$ for
all $z \in \mathbb{Z}$. 

Now consider the same group $G = \mathbb{Z}$ but this time with choice of generators $Y = \lbrace a,a',b \rbrace$ where $a$ and $a'$ both represent $+1$ and $b$ again represents $-1$.
Just as before $\lambda^{G,Y}_0(n) = 0$ for $n$ odd, and a word of length $n = 2k$ represents $0$ if and only if it contains $b$ in exactly $k$ positions. There are again $\binom{2k}{k}$ ways to choose the positions of the $b$s, but this time there are $2^k$ ways to choose whether the remaining $k$ positions should contain $a$ or $a'$. So using Stirling's approximation again:
$$\lambda_0^{G,Y}(n) \ = \ \binom{2k}{k}  2^k \ \approx \ \frac{4^k}{\sqrt{\pi k}} 2^k \ = \ \frac{1}{\sqrt{\pi \frac{n}{2}}} (2 \sqrt{2})^{n}$$
which gives $\lambda_0^{G,Y} = 2 \sqrt{2}$, which is not maximal. Using Proposition~\ref{prop_jorder} again, $\lambda_z^{G,Y} = 2 \sqrt{2} < 3$ 
all $z \in \mathbb{Z}$ so there are no elements of maximal local cogrowth.
\end{example}

\section{Global Cogrowth}

We saw above that, while local cogrowth displays interesting properties in some semigroups,
in other important semigroups it is everywhere $0$. This motivates another way to define cogrowth for
a semigroup, which does not involve looking at individual elements separately and which is non-trivial in these cases.

\begin{definition}
The \textit{global cogrowth function} of $S$ (with respect to the
choice of generators $X$)
is the function
$\gamma^{S,X} : \mathbb{N} \to \mathbb{N}$ which maps $n$ to the number
of pairs $(u,v)$ of non-empty words over $X$ with $|uv| = n$ and
$\ol{u} = \ol{v}$ in $S$.

The \textit{global cogrowth rate} $\gamma^{S,X}$ of $S$ is
$$\limsup_{n \to \infty} \ (\gamma^{S,X}(n))^{1/n}.$$

Again, we omit superscripts where no confusion can arise.
\end{definition}

As we start to study the global cogrowth rate in more depth, we shall see that it admits numerous equivalent formulations (see Remark~\ref{remark_localvsglobal}, Lemma~\ref{lemma_gammaprime}, Proposition~\ref{prop_sup} and Corollary~\ref{cor_cogrowthrandomwalk} below).

\begin{remark}
Clearly we have $0 \leq \gamma^{S,X}(n) \leq (n-1) |X|^n$ for all $n$. Hence, 
$0 \leq \left( \gamma^{S,X}(n) \right)^{1/n} \leq (n-1)^{1/n} |X|$ for all $n$, from which it follows
that $\gamma^{S,X}$ is defined and $0 \leq \gamma^{S,X} \leq |X|$.
\end{remark}

We will be interested not so much in the actual value of the global cogrowth rate, as in when it obtains its minimum and (especially) maximum
possible values.

\begin{definition}
We say that $S$ has \textit{maximal global cogrowth (with respect to $X$)} if $\gamma^{S,X} = |X|$.
\end{definition}

\begin{remark}\label{remark_localvsglobal}
As an immediate consequence of the definitions, we can describe global cogrowth in terms of local cogrowth by the formula:
$$\gamma^{S,X}(n) \ = \ \sum_{\substack{i, j \in \mathbb{N} \\ i+j = n \\ s \in S}} \lambda^{S,X}_s(i) \lambda^{S,X}_s(j).$$
Notice that the sum on the right hand side has only finitely many non-zero terms, since only finitely many elements can be
represented by words over $X$ of length $n$ or less.
\end{remark}

There is another elementary relationship between local and global cogrowth rates:

\begin{proposition}\label{prop_localvsglobal}
The global cogrowth rate is bounded below by the local cogrowth rate
of every element. In particular, every semigroup of maximal local cogrowth has maximal global cogrowth.
\end{proposition}
\begin{proof}
Let $s \in S$. For any $n$ there are $\lambda_s(n)$ words of length $n$
representing $s$ and hence at least
$\lambda_s(n)^2$ pairs $(u,v)$ of words both of which represent $s$ and
with $|uv| = 2n$. Thus, we have $\gamma(2n) \geq \lambda_s(n)^2$ for all $n$.
The claim follows.
\end{proof}

\begin{corollary}\label{cor_finitekesten}
Every finitely generated semigroup with a finite minimal left or right
ideal has maximal global cogrowth. In particular, every finite semigroup has
maximal global cogrowth.
\end{corollary}
\begin{proof}
This is immediate from Propositions~\ref{prop_minimalideal} and \ref{prop_localvsglobal}.
\end{proof}

\begin{example}
We saw in Example~\ref{example_bicycliclocal} that the ($2$-generated) bicyclic monoid has 
maximal local cogrowth at the identity element, so by Proposition~\ref{prop_localvsglobal}
it also has maximal global cogrowth.
\end{example}

\begin{proposition}\label{prop_globalquotient}
If $T$ is a homomorphic image of $S$ (with the obvious choice of generators obtained by composing the
map from $X^+$ to $S$ with the morphism from $S$ onto $T$) then the global cogrowth of $T$
is greater than or equal to the global cogrowth of $S$. In particular, the class of semigroups having maximal global
cogrowth with respect to some finite choice of generators is closed under the taking of quotients.
\end{proposition}
\begin{proof}
This follows immediately from the fact that two words representing
the same element of $S$ also represent the same element of $T$.
\end{proof}

\begin{example}\label{example_freecommutativeglobal}
Let $S$ be the free commutative monoid $\langle a, b \mid ab = ba \rangle$.
For any $n \in \mathbb{N}$ there are $\binom{2n}{n}$ distinct words of length $2n$ representing
$a^n b^n$, so using Stirling's approximation we obtain
$$\gamma(4n) \ \geq \ {\binom{2n}{n}}^2 \ \approx \ \left( \frac{4^n}{2n+1} \right)^2 \ = \ \frac{2^{4n}}{(2n+1)^2}$$
from which it follows (since $(2n+1)^2$ is subexponential) that $\gamma \geq 2$. Since there are only two generators
we must have $\gamma = 2$, that is, $S$ has maximal global cogrowth.
\end{example}
A comparison of Examples~\ref{example_freecommutativelocal} and \ref{example_freecommutativeglobal} illustrates
how greatly global cogrowth can differ from local cogrowth. In particular, a semigroup can have maximal global
cogrowth but local cogrowth which is everywhere $0$.

A similar, if slightly more involved, combinatorial analysis to that in
Example~\ref{example_freecommutativeglobal} can be used to compute the global
cogrowth of free commutative semigroups of all finite ranks; it turns out that they all have
maximal global cogrowth. Since every finitely generated commutative semigroup is a quotient of such
a semigroup, it follows by Proposition~\ref{prop_globalquotient} that
every finitely generated commutative semigroup has maximal global cogrowth. But
in fact we can prove something much stronger than this:

\begin{theorem}\label{thm_growthglobal}
Every finitely generated semigroup of subexponential growth has maximal global cogrowth.
\end{theorem}
\begin{proof}
Let $f : \mathbb{N} \to \mathbb{N}$ be the (spherical) growth function of
$S$ with respect to $X$, that is, the function which maps $n$ to the number
of distinct elements of $S$ which admit representatives in $X^+$ of length exactly $n$.
Let $k = |X|$.

Then for each $n$, by the pigeon hole principle, we may choose an element
$x_n \in S$ such that at least $\frac{k^n}{f(n)}$ distinct words of length
exactly $n$ represent $x_n$. Thus, there are at least
$\left( \frac{k^{   }}{f(n)} \right)^2$
ways to choose words $u$ and $v$ of length $n$ so that $u = v$ in $S$, and so
$$\gamma(2n) \ \geq \ \left( \frac{k^n}{f(n)} \right)^2 \ = \ \frac{1}{f(n)^2} k^{2n}$$
giving
\begin{align*}
\gamma \ &= \ \limsup_{n \to \infty} \ (\gamma(n))^{1/n} \\
&\geq \ \limsup_{m \to \infty} \ (\gamma(2m))^{1/2m}  \\
&\geq \ \limsup_{m \to \infty} \ \left( \frac{1}{f(m)^2} k^{2m} \right)^{1/2m} \\
&= \ \limsup_{m \to \infty} \ \left( \frac{1}{f(m)} \right)^{1/2m} k \\
&= \ k \ \limsup_{m \to \infty} \ \left( \frac{1}{f(m)} \right)^{1/2m} \\
&= \ k
\end{align*}
where the last equality holds because the growth function $f$ is subexponential.
\end{proof}

In particular, we obtain by this route the claimed fact about commutative semigroups:
\begin{corollary}\label{cor_commutativeglobal}
Every finitely generated commutative semigroup has maximal global cogrowth.
\end{corollary}

\begin{remark}
The converse to Theorem~\ref{thm_growthglobal} fails even for groups with symmetric choices of generators; indeed by Theorem~\ref{thm_groupamen} below, any finitely generated  amenable group of exponential growth (for example, any solvable group which is not virtually nilpotent, such as the Baumslag-Solitar group $BS(1,2)$) with any finite symmetric generating set provides a counterexample.
\end{remark}

\begin{lemma}\label{lemma_monoidideal}
Let $S$ be a semigroup with choice of generators $X$, having an ideal $I$ which is a monoid with identity element $e$. Consider $X$ as a
choice of generators for $I$, with each symbol $x \in X$ representing the element $e\ol{x} \in I$. Then
$\gamma^{S,X} = \gamma^{I,X}$.
\end{lemma}
\begin{proof}
The map
$$S \to I, \ x \mapsto ex$$
is easily verified to be a surjective homomorphism from $S$ onto $I$, so by
Proposition~\ref{prop_globalquotient} we have $\gamma^{I,X} \geq \gamma^{S,X}$.

Conversely, fix a word $w$ over $X$ representing $e$ in $S$. Then
for any words $u$ and $v$, we have $u \equiv v$ in $I$ if and only if
$wu \equiv wv$ in $S$. Hence for every $n \in \mathbb{N}$ we have
$$\gamma^{I,X}(n) \ \leq \ \gamma^{S,X}(n+2|w|)$$
and so $\gamma^{S,X} \leq \gamma^{I,X}$.
\end{proof}

We now introduce a slight variation of the definition of the global cogrowth function, which will turn out to be
equivalent for the purpose of calculating the cogrowth rate, but is technically easier to work with since it
involves comparing only words of the same length:

\begin{definition} Let $S$ be a semigroup with finite choice of generators $X$. We define
$\gamma' : \mathbb{N} \to \mathbb{N}$ by 
$$\gamma'(n) \ = \ \left| \lbrace (u,v) \mid u, v \in X^+, u \equiv v, |u| = |v|, |uv| = n \rbrace \right|$$
\end{definition}

\begin{remark}
Immediately from the definition, $\gamma'(n) = 0$ for $n$ odd, while for $n = 2k$ we have
$$\gamma'(n) = \sum_{s \in S} \lambda_s(k)^2.$$
\end{remark}

\begin{lemma}\label{lemma_gammaprime}
Let $S$ be a semigroup finitely generated by $X$. Then
$$\gamma \ = \ \limsup_{n \to \infty} \left( \gamma'(n) \right)^{1/n} \ = \ \limsup_{\substack{n \to \infty \\ n \textrm{ even}}} \left( \gamma'(n) \right)^{1/n}$$
\end{lemma}

\begin{proof}
The inequalities
$$\gamma \ \geq \ \limsup_{n \to \infty} \left( \gamma'(n) \right)^{1/n} \ \geq \ \limsup_{\substack{n \to \infty \\ n \textrm{ even}}} \left( \gamma'(n) \right)^{1/n}$$
are clear. Indeed, the first holds because by definition $\gamma(n) \geq \gamma'(n)$ for all $n$, and the second from the definition of the limit superior.

For the remaining inequality, define
$$\gamma' \ = \ \limsup_{\substack{n \to \infty \\ n \textrm{ even}}}  \left( \gamma'(n) \right)^{1/n}$$
and suppose for a contradiction that $\gamma' < \gamma$. Choose
some $\delta$ with $\gamma' < \delta < \gamma$. Then there is a constant $\beta$ such that
$\gamma'(n) \leq \beta \delta^{n}$ for all even $n$. From the definition of $\gamma'$
this means that for every $i$ we have
$$\sum_{s \in S} \lambda_s(i)^2 \ = \ \gamma'(2i) \ \leq \ \beta \delta^{2i}.$$
Recall from Remark~\ref{remark_localvsglobal} that
$$\gamma(n) \ = \ \sum_{\substack{i+j = n \\ s \in S}} \lambda_s(i) \lambda_s(j) \ = 
\ \sum_{i+j = n} \left( \sum_{s \in S} \lambda_s(i) \lambda_s(j) \right).$$
Notice that, since for each $k$ there are only finitely many $s \in S$ with
$\lambda_s(k) \neq 0$, each term in the outer sum is the scalar product
of two finite dimensional vectors; by the Cauchy-Schwartz inequality this cannot exceed the product of their Euclidean norms, so whenever $i+j=n$ we have:
$$\sum_{s \in S} \lambda_s(i) \lambda_s(j) \ \leq \  \sqrt{\sum_{s \in S} \lambda_s(i)^2} \sqrt{\sum_{s \in S} \lambda_s(j)^2} \ \leq \ \sqrt{\beta \delta^{2i}} \sqrt{\beta \delta^{2j}} \ = \ \beta \delta^n$$
and hence, since there are $n-1$ possible choices of $i$ and $j$,
$$\gamma(n) \ \leq \ \sum_{i+j = n} \beta \delta^n \ \leq \ (n-1) \beta \delta^{n}.$$
It follows that $\gamma \leq \delta$, contradicting the choice of $\delta$
and completing the proof.
\end{proof}

\begin{remark}\label{remark_globalrandom}
The definition of $\gamma'(n)$ allows an interpretation of global cogrowth in terms of random walks, similar to
that for local cogrowth described in Remark~\ref{remark_localrandomwalk}.  It is easy to see that $\gamma'(2n) / |X|^{2n}$ is the probability
that two independent random walks, each starting at the identity element (or at an adjoined identity element for
semigroups without identity) and proceeding by multiplying consecutively on the right (say) by $n$ generators chosen independently and uniformly
at random, end at the same element. By Lemma~\ref{lemma_gammaprime}, the global cogrowth is maximal if and only if this
probability does not decay exponentially fast as $n \to \infty$.
\end{remark}

In the special case of a group equipped with a symmetric choice of (semigroup) generators, it transpires that there is no essential difference between local and global cogrowth:
\begin{proposition}\label{prop_localglobalgroup}
Let $G$ be a group with a finite symmetric choice of generators $X$. Then
$$\gamma^{G,X} = \lambda^{G,X}_g$$
for all $g \in G$.
\end{proposition}
\begin{proof}
Let $1$ denote the identity element of $G$. We shall need a notion of a formal inverse of a word over the generating set $X$. Since the choice of generators is symmetric, we may pair up each generator $x \in X$ with a generator $x' \in X$ in such a way that $x$ and $x'$ represent mutually inverse elements and $(x')' = x$. (If $x$ represents the identity or an involution it may be necessary to choose $x' = x$.) Extend to words in the obvious way, by defining $(x_1 \dots x_n)' = (x_n)' \dots (x_1)'$ so that $u$ and $u'$ always represent
mutually inverse elements.

By Proposition~\ref{prop_jorder} we have $\lambda_g = \lambda_h$ for all $g, h \in G$, so it 
suffices to show that $\lambda_1 = \gamma$. In fact, by Lemma~\ref{lemma_gammaprime}, it suffices to show that
$\lambda_1(2n) = \gamma'(2n)$ for all $n$.

But every word of length $2n$ can be written uniquely in the form $u (v')$ for some 
$u, v \in X^+$ with $|u| = |v| = n$, and clearly we have $u \equiv v$ in $G$ if and only if $u (v') \equiv 1$
in $G$. Thus, there is a one-one correspondence between the pairs counted by $\gamma'(2n)$ and the
words counted by $\lambda_1(2n)$.
\end{proof}

The hypothesis in Proposition~\ref{prop_localglobalgroup} that the choice of generators be symmetric cannot be
removed, as the following example shows:
\begin{example}\label{example_groupnonsymmetric}
Recall from Example~\ref{example_groupnonsymmetric_local} that when $G = \mathbb{Z}$ with two generators representing $+1$ and only one representing
$-1$ we have $\lambda_z = 2 \sqrt{2}$ for all $z \in \mathbb{Z}$. On the other hand, Corollary~\ref{cor_commutativeglobal} tells us that $\mathbb{Z}$ has maximal global cogrowth
irrespective of the semigroup generating set chosen, so for this generating set $\gamma = 3$.
\end{example}

We next consider the case of (non-commutative) free semigroups; it turns out that these are (except for a pathological case where $|X|=1$) exactly characterized by their global cogrowth:

\begin{theorem}\label{thm_free}
Let $S$ be a semigroup with finite choice of generators $X$. Then
$$\gamma^{S,X} \ \geq \ \sqrt{|X|}$$
with equality if and only if either $S$ is a free
semigroup freely generated by $X$ or $|X| = 1$.
\end{theorem}

\begin{proof}
First suppose $S$ is a free semigroup, freely generated by $X$.
For any words $u, v \in X^+$ we
have $\ol{u} = \ol{v}$ in $S$ if and only
if $u=v$ as words. It follows easily that
$$\gamma(n) \ = \ \begin{cases}
0 &\textrm{ if $n$ is odd;} \\
|X|^{n/2} &\textrm{ if $n$ is even}
\end{cases}$$
whence
$$\gamma \ = \ \limsup_{n \to \infty} \gamma(n)^{1/n}
\ = \ \limsup_{m \to \infty} \gamma(2m)^{1/2m}
\ = \ \limsup_{m \to \infty} (|X|^m)^{1/2m}
\ = \ \sqrt{|X|}.$$

Now consider the case where $S$ is not (necessarily) free on $X$.
Since every $k$-generated semigroup is a quotient of a free semigroup of rank $k$,
it follows from this and Proposition~\ref{prop_globalquotient} that
$\sqrt{k}$ is the minimum possible value for the global cogrowth of
a $k$-generated semigroup.

In particular, if $|X| = 1$ then we have $\gamma \geq \sqrt{|X|} = 1$, and clearly $\gamma \leq |X| = 1$ so in this case $\gamma = 1 = |X|$.

Finally, suppose $|X| \geq 2$ and $S$ is not freely generated by $X$. Then there are two distinct words $u, v \in X^+$ which represent the same element of $S$. We
claim we may choose these words to have the same length. Indeed, if $u$ and $v$ do not commute then it suffices to replace $u$ by $uv$ and $v$ by
$vu$. If they do commute, then by \cite[Proposition 1.3.2]{lothaire_algebraic} they are powers of a common word, say $u = w^q$ and $v = w^r$. Now let $a$ and $b$ be
distinct symbols in $X$ (recalling that $|X| \geq 2$) and 
choose $n$ larger than the lengths of $u$ and $v$. Let $u' = w^q ab^{n+2}$ and $v' = w^r ab^{n+2}$. Clearly $u'$ and $v'$ are distinct words representing the same
element of $S$. It is easily seen that neither $u'$ nor $v'$ can be a proper power: indeed, if either was then it would have to be a power of a suffix of $b^{n+2}$ (since the latter
makes up more than half the word) but this contradicts the fact it contains the letter $a$. Thus, they cannot commute, so replacing $u$ with $u'v'$ and $v$ with $v'u'$ yields the required properties. 
Let $p$ be the common length of $u$ and $v$.

Consider now a pair $(a,b)$ of two (not necessarily distinct) words $a, b \in X^+$ of length $n$, both representing the same element of $S$. From $(a,b)$ we can construct two different types of pairs of
words of length $n+p$:
\begin{itemize}
\item the $|X|^p$ pairs of the form $(aw,bw)$ where $w \in X^+$ is of length $p$;
\item the $2$ pairs $(au,bv)$ and $(av,bu)$.
\end{itemize}
Clearly, each pair consists of two words representing the same element of $S$.
The pairs of the second kind are distinct from those of the first kind since the two sides do not share a suffix of length $p$, while the pairs within each kind are distinct
from each other because the final $p$ letters of (say) the left-hand-side are always different. Thus, we have constructed $|X|^p+2$ distinct pairs. Moreover, the words $a$ and $b$
are recoverable as the $n$-letter prefixes of all the pairs, so a different choice of $(a,b)$ of the same length would clearly lead to a disjoint collection of pairs. Thus, we have
constructed $|X|^p+2$ pairs of words of length $n+p$ for each pair of words of length $n$, which means that
$$\gamma'(2(n+p)) \ \geq \ (|X|^p + 2) \ \gamma'(2n)$$
for all $n$. A simple induction gives
$$\gamma'(2+2kp) \ \geq \ (|X|^p + 2)^k \gamma'(2) \ \geq \ (|X|^p + 2)^k$$
for all $k$, and it follows easily by Lemma~\ref{lemma_gammaprime} that
$$\gamma \ \geq \ (|X|^p + 2)^{1/2p} \ > \ (|X|^p)^{1/2p} \ = \ \sqrt{|X|}.$$
\end{proof}

We now return to studying the properties of the global cogrowth functions $\gamma$ and $\gamma'$. The following two lemmas give a sense in which these functions approach their limiting growth rate in a reasonably uniform way; this will be of great importance for many of our subsequent results.

\begin{lemma}\label{lemma_globalcogrowthsmooth1}
If $\gamma'(n_0) > \kappa^{n_0}$ for some $\kappa \geq 0 $ and $n_0 \in \mathbb{N}$ then there is a constant $C > 0$ such that $\gamma'(n) > C \kappa^n$ for all
even $n \in \mathbb{N}$.
\end{lemma}
\begin{proof}
If $\kappa = 0$ then the claim is trivial, so suppose $\kappa > 0$.
Suppose $\gamma'(n_0)  > \kappa^{n_0}$. Since $\gamma'(n) = 0$ for odd $n$, we must have $n_0$ even, say $n_0 = 2q$.

First, we claim that $\gamma'(2pq) > \kappa^{2pq}$ for all $p \in \mathbb{N}$. Indeed, by the definition of $\gamma'$ there are 
$\gamma'(2q)$ pairs of words of length $q$ representing the same element, so by choosing $p$ such pairs and concatenating
the respective sides we can obtain $\left( \gamma'(2q) \right)^p > \kappa^{2pq}$ pairs of words of length $pq$ representing the same
element. Thus, $\gamma'(2pq) > \kappa^{2pq}$ as claimed.

Now for any $p \geq 1$ and $1 \leq r < q$ there are at least as many pairs of words of length $pq+r$ representing the same element as
there are pairs of words of length $pq$ representing the same element.  Indeed, from each pair of the latter kind one may obtain a pair of
the former kind by fixing a letter $a$ and concatenating $a^r$ to each side. Thus, we have 
$$\gamma'(2(pq+r)) \ \geq \ \gamma'(2pq) \ > \ \kappa^{2pq} \ = \ \kappa^{-2r} \kappa^{2(pq+r)}.$$
So if we choose $C$ smaller than the finitely many values of $\kappa^{-2r}$ for $1 \leq r < q$ and the finitely many values of $\gamma'(2m)  \kappa^{-2m}$ for
$1 \leq m < q$, we have $\gamma'(2m) > C \kappa^{2m}$ for all $m \in \mathbb{N}$, as required.
\end{proof}

\begin{lemma}\label{lemma_globalcogrowthsmooth2}
For any $0 \leq \kappa < \gamma$, there is a constant $C > 0$ such that
$$\gamma(n) \ \geq \ \gamma'(n) \ > \ C \kappa^n$$
for all even $n \in \mathbb{N}$.
\end{lemma}
\begin{proof}
The first inequality is immediate from the definitions.

For the second, 
since $\kappa < \gamma$, it follows by Lemma~\ref{lemma_gammaprime} that there is an $n \in \mathbb{N}$ with $\gamma'(n) > \kappa^{n}$,
and then the claim is immediate from Lemma~\ref{lemma_globalcogrowthsmooth1}.
\end{proof}

A consequence of Lemma~\ref{lemma_globalcogrowthsmooth1} is that we have the option of characterising global cogrowth as a simple supremum instead of a limit superior:
\begin{proposition}\label{prop_sup}
Let $S$ be a semigroup with finite choice of generators $X$. Then 
$$\gamma \ = \ \sup_{n \in \mathbb{N}} \left( \gamma'(n) \right)^{1/n} \ = \ \sup_{\substack{n \in \mathbb{N} \\ n \textrm{ even}}} \left( \gamma'(n) \right)^{1/n}.$$
\end{proposition}
\begin{proof}
We prove first the left-hand equality. By Lemma~\ref{lemma_gammaprime} we have
$$\gamma  \ = \ \limsup_{n \to \infty} \left( \gamma'(n) \right)^{1/n}.$$
and by definition the limit superior cannot exceed the supremum.

Now  suppose for a contradiction that the supremum in the middle strictly
exceeds $\gamma$, and choose $\kappa$ with 
$$\gamma \ < \ \kappa \ < \ \sup_{n \in \mathbb{N}} \left( \gamma'(n) \right)^{1/n}.$$
Then there is an $n_0 \in \mathbb{N}$ such that  $\kappa < \left( \gamma'(n_0) \right)^{1/n_0}$,
that is $\kappa^{n_0} < \gamma'(n_0)$. Now by Lemma~\ref{lemma_globalcogrowthsmooth1}
there is a constant $C > 0$ so that $\gamma'(n) > C \kappa^n$ for all even $n > 1$, from which
it follows that $\gamma \geq \kappa$, giving the required contradiction. 

The right-hand equality in the statement is immediate from the fact that $\gamma'(n)$ is $0$ for odd $n$ and non-negative for even $n$.
\end{proof}

Our first application of the above uniformity results is to show that, where a semigroup has maximal cogrowth with respect to
some choice of generators, one is free to assume (passing to the monoid $S^1$ if necessary) that the choice of generators contains
a representative for the identity element. (Recall that, by definition, if $S$ is a monoid then $S = S^1$.)

\begin{theorem}\label{thm_addidentity}
Suppose $S$ is a semigroup or monoid with finite choice of generators $X$, and let $Y$ be the choice of generators for the monoid $S^1$ obtained by adding an extra generator representing the identity element of $S^1$. Then $\gamma^{S^1,Y} \geq \gamma^{S,X} + 1$. In particular, if $S$ has maximal global cogrowth with respect to
$X$ then $S^1$ has maximal global cogrowth with respect to $Y$.
\end{theorem}
\begin{proof}
First note that if $|X| = 1$ then $S$ has subexponential growth, so $S^1$ has subexponential growth, which by Theorem~\ref{thm_growthglobal} means it has maximal cogrowth with respect to every generating set, so that $\gamma^{S^1,Y} = |Y| = |X|+1 = \gamma^{S,X} + 1$.
We may assume, therefore that $|X| > 1$. It follows by Theorem~\ref{thm_free} that $\gamma^{S,X} \geq \sqrt{|X|} > 1$.

Let $e$ be a new symbol representing the identity element of $S^1$, and consider the choice of generators $Y = X \cup \lbrace e \rbrace$ for $S^1$.

Let $n \in \mathbb{N}$, and consider the value of $\gamma^{S^1,Y}(n)$, that is, the number of pairs of words of total length $n$ over $Y$ representing
the same element in $S^1$. Clearly, since inserting $e$s into a word does not change the element of $S^1$ represented, each such pair may be
obtained from a unique (typically shorter) pair over the alphabet $X$ by inserting $e$s at different points in the words to increase their length to $n$.
Given a pair of words in $X^+$ of total length $n-i$, there are at least\footnote{In fact there are $n-i+2$ but using only $n-i+1$ of them gives a cleaner argument.} $n-i+1$ possible positions into which to insert the required $i$ $e$'s,
so the number of ways to do this is at least
$$\binom{i + (n-i+1) - 1}{(n-i+1)-1} = \binom{n}{n-i} = \binom{n}{i}.$$
Thus, every pair of words of total
length $n-i$ over $X$ yields $\binom{n}{i}$ pairs of words of total length $n$ over $Y$, so we have 
$$\gamma^{S^1,Y}(n) \ \geq \ \sum_{i=0}^n \binom{n}{i} \gamma^{S,X}(i).$$
Recalling that $\gamma^{S,X} > 1$, let $1 \leq \kappa < \gamma^{S,X}$. By Lemma~\ref{lemma_globalcogrowthsmooth2}, there is a $C > 0$ such that
$\gamma^{S,X}(i) > C \kappa^i$ for all
even $i$. Thus, we have
$$\gamma^{S^1,Y}(n) \ \geq \ \sum_{i=0}^n \binom{n}{i} \gamma^{S,X}(i) \ \geq \ \sum_{\substack{i = 0 \\ i \textrm{ even}}}^n \binom{n}{i} C \kappa^i \ = \ C \sum_{\substack{i = 0 \\ i \textrm{ even}}}^n \binom{n}{i} \kappa^i.$$

Now, for convenience defining $\binom{n-1}{-1} = \binom{n-1}{n} = 0$, we have
\begin{align*}
\sum_{\substack{i = 0 \\ i \textrm{ even}}}^n \binom{n}{i} \kappa^i
\ &= \ \sum_{\substack{i = 0 \\ i \textrm{ even}}}^n \left( \binom{n-1}{i-1} + \binom{n-1}{i} \right) \kappa^i & \\
&= \ \left( \sum_{\substack{i = 0 \\ i \textrm{ even}}}^n \binom{n-1}{i} \kappa^i \right) + \left( \sum_{\substack{j = -1 \\ j \textrm{ odd}}}^{n-1} \binom{n-1}{j} \kappa^{j+1} \right) \\
&\geq \ \left( \sum_{\substack{i = 0 \\ i \textrm{ even}}}^n \binom{n-1}{i} \kappa^i \right) + \left( \sum_{\substack{j = -1 \\ j \textrm{ odd}}}^{n-1} \binom{n-1}{j} \kappa^{j} \right) \\
&= \ \sum_{i = 0}^{n-1} \binom{n-1}{i} \kappa^{i} \\
&= \ (\kappa+1)^{n-1} 
\end{align*}
where the inequality holds because $\kappa > 1$ so $\kappa^{j} < \kappa^{j+1}$, and the final equality is an application of the binomial theorem.
Thus for all $n$ we have
$$\gamma^{S^1,Y}(n) \geq C (\kappa+1)^{n-1} \ = \ \left( C (\kappa+1)^{-1} \right)  (\kappa+1)^n.$$
Hence $\gamma^{S^1,Y} > \kappa+1$, and since $\kappa$ was an arbitrary value less than $\gamma^{S,X}$ we must have
$\gamma^{S^1,Y} \geq \gamma^{S,X}+1$ as required.
\end{proof}

\begin{theorem}\label{thm_supersets}
Suppose $S$ is a semigroup or monoid with maximal global cogrowth with respect to some choice of generators, and let $K$ be any finite
subset of $S^1$. Then the monoid $S^1$ has maximal cogrowth with respect to some choice of generators containing a representative
for every element of $K$.
\end{theorem}

\begin{proof}
By Theorem~\ref{thm_addidentity} there is a choice of generators $X$ for $S^1$ containing a generator, say $e$, representing the identity element. Since $K$ is finite and $X$ generates $S^1$, there is a $p$
such that every element of $K$ can be represented by a word over $X$ of length at most $p$. Moreover, since $X$ contains a generator representing
the identity, every element of $K$ can be represented by a word over $X$ of length exactly $p$, as can each element represented by a generator
from $X$.

Thus, there is a choice of generators (call it $Y$) for $S^1$ containing one generator for each word of length $p$ over $X$. Now for $n$ divisible by $2p$, say $n = 2pq$, every pair of words over $X$ of length $pq$ representing the same element of $S^1$ clearly yields a distinct pair of words over $Y$ of length $q$ representing the same element of $S^1$. Hence, we have
$$\gamma'^{S^1,Y}(2q) \ \geq \ \gamma'^{S^1,X}(2pq)$$
for all $q$.

Now for any $\kappa < |X| = \gamma^{S^1,X}$ by Lemma~\ref{lemma_globalcogrowthsmooth2} there is a $C > 0$ so that
$\gamma'^{S^1,X}(n) > C \kappa^n$ for all even $n$. Thus, for all even $n$, say $n = 2q$, we have
$$\gamma'^{S^1,Y}(n) \ = \ \gamma'^{S^1,Y}(2q) \ \geq \ \gamma'^{S^1,X}(2pq) \ > \ C \kappa^{2pq} \ = \ C (\kappa^{p})^{n}$$
so by Lemma~\ref{lemma_gammaprime} we have $\gamma^{S^1,Y} \geq \kappa^p$. Since $\kappa$ was an arbitrary value less than $|X|$ it follows
that $\gamma^{S^1,Y} \geq |X|^p = |Y|$ so $S^1$ has maximal cogrowth with respect to the choice of generators $Y$, as required.
\end{proof}

\section{Cogrowth and the Markov Operator of a Random Walk}\label{sec_randomwalk}

Consider again a random walk, starting at the identity of the monoid $S^1$ and multiplying on the right by generators chosen uniformly
at random from $X$.  This is a Markov process, and has an associated $S^1 \times S^1$ transition matrix $M$ where for $s,t \in S^1$ the entry $M_{s,t}$
--- the probability that starting at $s$ and taking one step takes one to $t$ --- is given by
$$M_{s,t}  = \frac{| \lbrace x \in X \mid s \ol{x} = t \rbrace |}{|X|}.$$
Being clearly row stochastic, this matrix is an operator on the right of the Banach space $\ell_1(S^1)$.

(Of course we could also, dually, consider a random walk obtained by multiplying on the \textit{left} by generators, obtaining a Markov operator
on the left of $\ell_1(S^1)$. In this section we prefer to work with the right random walk operator, for consistency with antecedent work of Day \cite{Day64} and with the bulk of the literature on semigroup Cayley graphs. Of course, all of our results admit left-right duals, the statements of which are easily obtained and some of which we shall need to use in Section~\ref{sec_amenable} below.)

Let $v \in \ell_1(S^1)$ be the probability mass function corresponding to some probability distribution on $S^1$. Then $vM^n$ is the probability
mass function resulting from starting in the distribution corresponding to $v$ and taking $n$ steps of the random walk. Notice also that
the $\ell_2$-norm
$$|v|_2 = \sqrt{\sum_{s \in S^1} v(s)^2}$$
is the square root of the probability that two elements of $S^1$, selected independently at random according to $v$, coincide. If we let $\chi_1 \in \ell_1(S^1)$ be
the characteristic function of the identity element in $S^1$ (that is, the probability mass function of the starting distribution for our random walk)
then these observations combine with Remark~\ref{remark_globalrandom} to prove the following connection between the operator $M$ and
global cogrowth:
\begin{proposition}\label{prop_cogrowthrandomwalk}
For all $n$, $\gamma'(2n) \ = \ |X|^{2n} \left( |\chi_1 M^n|_2 \right)^2$.
\end{proposition}
Combining with Lemma~\ref{lemma_gammaprime} and Proposition~\ref{prop_sup} this yields:
\begin{corollary}\label{cor_cogrowthrandomwalk}
$$\gamma \ = \ |X| \limsup_{n \to \infty} \left( |\chi_1 M^n|_2 \right)^{1/n} \ = \ |X| \sup_{n \in \mathbb{N}} \left( |\chi_1 M^n|_2 \right)^{1/n}.$$
\end{corollary}

The natural way in which the $\ell_2$-norm arises suggests that it might be profitable to study $M$ as an operator on the Hilbert space $\ell_2(S^1)$. In general, unfortunately, $M$ is not defined on the whole of $\ell_2(S^1)$, and even when it does define a map from $\ell_2(S^1)$ to functions $S^1 \to \mathbb{R}$ it may not preserve square-summability. We shall shortly (Proposition~\ref{prop_operatorcondition} below) describe exactly the circumstances under which it does act on $\ell_2(S^1)$,
but first we consider a special case.

Suppose for now that $M$ is an operator on $\ell_2(S^1)$, and recall that the \textit{spectral radius} of $M$ is the supremum of the absolute values of all $\lambda$ such that $M - \lambda I$ is singular, where
$I$ is the identity operator.
By the well-known theorem of Gelfand, provided $M$ is bounded this value is equal to the rate of exponential growth of
the operator norm of the powers of $M$, that is,
$$\limsup_{n \to \infty} \left( ||M^n||_2 \right)^{1/n},$$
where $||N||_2$ denotes the operator norm
$\sup_{v \in \ell_2(S^1) \setminus \lbrace 0 \rbrace} \frac{|vN|_2}{|v|_2}$ of an operator $N$ on $\ell_2(S^1)$.

In the case $S$ is a group (hence $S = S^1$) and $X$ is a symmetric choice of generators, a celebrated theorem of Kesten \cite{Kesten59b} describes
\textit{amenability} of the group $S$: 
\begin{theorem}[Kesten 1959]\label{thm_kesten}
Let $S$ be a group, $X$ a symmetric choice of generators for $S$, and $M$ the corresponding right random walk transition matrix. Then
\begin{itemize}
\item $M$ is a bounded symmetric operator on $\ell_2(S)$ with norm (and hence spectral radius) at most $1$; and
\item $M$ has spectral radius $1$ if and only if $S$ is amenable.
\end{itemize}
\end{theorem}

We now return to a more general semigroup $S$ and the question of whether or when $M$ is an operator (or even a bounded operator) on $\ell_2(S^1)$.
Recall that a semigroup has \textit{bounded right indegree} (sometimes called \textit{finite geometric type}) if for every element $x \in S$ there is a $b \in \mathbb{N}$ such that for each
$s \in S$ no more than $b$ elements $t \in S$ satisfy $tx = s$. In a finitely generated semigroup, this is easily seen to mean exactly that there is 
a bound on the number of edges coming into each vertex of the right Cayley graph. We shall need the following elementary inequality which can be proved in several different ways from the Cauchy-Schwarz inequality:

\begin{lemma}\label{lemma_squaresum}
For any $s_1, \dots, s_n \in \mathbb{R}$,
$$(s_1 + \dots + s_n)^2 \leq n (s_1^2 + \dots + s_n^2).$$
\end{lemma}
\begin{proof}
For $1 \leq i \leq n$, define $v_i \in \mathbb{R}^n$ by $(v_i)_j = s_{(j-i+1)}$. By direct calculation,
$(s_1 + \dots + s_n)^2 = (v_1 \cdot v_1) + (v_1 \cdot v_2) + \dots + (v_1 \cdot v_n)$.
Now for each $i$ the Cauchy-Schwarz inequality gives, $v_1 \cdot v_i \leq v_1 \cdot v_1 = v_i \cdot v_i$, so this means
$(s_1 + \dots + s_n)^2 \leq n (v_1 \cdot v_1) = n (s_1^2 + \dots + s_n^2).$ as required.
\end{proof}

\begin{proposition}\label{prop_operatorcondition}
Let $S$ be a finitely generated semigroup and let $M$ be the associated right random walk transition matrix as defined above. Then the following are equivalent:
\begin{itemize}
\item[(i)] $M$ defines an operator on $\ell_2(S^1)$;
\item[(ii)] $M$ defines a bounded operator on $\ell_2(S^1)$;
\item[(iii)] $S$  has bounded right indegree.
\end{itemize}
Moreover, if $S$ is right cancellative then $||M||_2 \leq 1$, and hence $M$ has spectral radius at most $1$.
\end{proposition}

\begin{proof}
The implication (ii) $\implies$ (i) being immediate, it suffices to show (iii) $\implies$ (ii) and (i) $\implies$ (iii).

(iii) $\implies$ (ii). Suppose $S$ has bounded right indegree. Let $v \in \ell^2(S)$. We need to show that $vM$ is defined and square-summable. For
each generator $a \in X$, let $b_a$ be a corresponding bound on the right indegree, and let $M_a$ be the matrix with
$$(M_a)_{st} = \begin{cases} 1 &\textrm{ if } sa = t \\  0 &\textrm{ otherwise} \end{cases}$$
so that $M = \frac{1}{|X|} \sum_{a \in X} M_a$.

Now for any $s \in S^1$
$$(v M_a)(s) = \sum_{ta = s} v(t)$$
is clearly defined, since there are at most $b_a$ elements $t$ satisfying $ta = s$, so $v M_a$ is defined as a function $S^1 \to \mathbb{R}$.
Moreover, 
$$(|v M_a|_2)^2 \ = \ \sum_{s \in S} ((v M_a)(s))^2 \ = \ \sum_{s \in S} \left( \sum_{t \in s a^{-1}} v(t) \right)$$
where $s a^{-1}$ denotes the set $\lbrace t \in S^1 \mid ta = s \rbrace$. For each $s$, this set has at most $b_a$ elements. Moreover, since for a
fixed $a$ and $t$ there is only one $s$ satisfying $ta = s$, these sets are disjoint as $s$ varies.
Hence $(|v M_a|_2)^2$ can be written as a sum of terms of the form
$$|v(t_1) + \dots + v(t_c)|^2 \ \leq \ b_a |v(t_1)|^2 + \dots + b_a |v(t_c)|^2$$
where $c \leq b_a$ so that the inequality follows from Lemma~\ref{lemma_squaresum}, and each $t \in S$ occurs in at most one term. Thus, we have
$$(|v M_a|_2)^2 \ \leq \ \sum_{t \in S} b_a |v(t)|^2 \ \leq \ b_a (|v|_2)^2$$
so that $|v M_a|_2 \leq \sqrt{b_a} |v|_2$. Thus, $M_a$ is a bounded operator on $\ell_2(S^1)$ and $||M_a||_2 \leq \sqrt{b_a}$. Hence by the triangle inequality:
$$||M||_2 \ = \ ||\frac{1}{|X|} \sum_{a \in X} M_a ||_2 \ \leq \ \frac{1}{|X|} \sum_{a \in X} \sqrt{b_a}.$$
so $M$ is a bounded operator on $\ell_2(S^1)$.

(i) $\implies$ (iii). We prove the contrapositive. Suppose that $S$ has unbounded indegree; since the generating set $X$ is finite, this means we may choose an $x \in X$
and a sequence of (not necessarily distinct) elements $s_1, s_2, \dots \in S$ such that for each $i$, there are distinct elements $t_{i,1}, \dots, t_{i,i}$ with $t_{i,j} x = s_i$ for all $j$.
By skipping some terms in the sequence of $s_i$'s if necessary, we may choose things so that the $t_{i,j}$'s are all distinct (as both $i$ and $j$ vary). We shall define an element $v \in \ell_2(S)$
such that $vM$ is not square-summable.

Define $v : S^1 \to \mathbb{R}$ by $v(t_{i,j}) = i^{-1.1}$ and $v(t) = 0$ for $t$ not of the form $t_{i,j}$. Then
$$\sum_{s \in S} |v(s)|^2 \ = \ \sum_{j \leq i} (t_{i,j})^2 \ = \ \sum_{i} i \left( i^{-1.1} \right)^2 \ = \ \sum_{i} i^{-1.2}$$
converges, so $v \in \ell_2(S^1)$. Notice that, because $v$ is supported only on the elements $t_{i,j}$ and because $t_{i,j} x = s_i$ for all $i$ and $j$,
we have that $v M_x$ is supported only on the elements $s_i$, and
$$(v M_x)(s_i) \ = \ \sum_{j \leq i} v(t_{i,j}) \ = \ i (i^{-1.1}) \ = \ i^{-0.1}.$$
Thus,
$$\sum_{s \in S} |(vM_x)(s)|^2 \ = \ \sum_{i} |(vM_x)(s_i)|^2 \ = \ \sum_i (i^{-0.1})^2 \ = \ \sum_i i^{-0.2},$$
which diverges, so $v M_x \notin \ell_2(S^1)$. It follows easily that $v M \notin \ell_2(S^1)$, so $M$ is not
an operator on $\ell_2(S^1)$, completing the proof that $(i) \implies (ii)$.

For the final part of the statement, notice that if $M$ is right cancellative then in the argument for (iii) $\implies$ (ii) we may take $b_a = 1$ for all $a$, giving
$$||M||_2 \ \leq \ \frac{1}{|X|} \sum_{a \in X} \sqrt{1} \ = \ 1.$$
It is well known and follows easily from the definitions that the spectral radius is bounded above by the operator norm. 
\end{proof}

\begin{remark}\label{remark_normintuition}
Intuitively, the norm and spectral radius of the Markov operator measure the rate at which the right random walk diffuses around the semigroup,
in one step (norm) or asymptotically (spectral radius). Values less than $1$ correspond to rapid diffusion, while values greater than $1$ imply rapid
concentration (which, as is made precise in the last part of Proposition~\ref{prop_operatorcondition}, can only result from a failure of right cancellativity). A value of exactly $1$ means that any diffusion or concentration is ``slow''. Kesten's Theorem therefore says, intuitively, that amenable groups are those in which random walks do not (asymptotically) diffuse rapidly across the group. Cases where the operator is undefined result from ``unboundedly fast concentration'', and the undefined operator can be thought of as having ``infinite norm and spectral radius''. For this reason, such cases are usually naturally considered together with cases where the operator is defined and has large (greater than $1$) norm and spectral radius. Many of our results will therefore feature conditions which are conjunctions such as ``$M$ is undefined or has spectral radius at least $1$''.
\end{remark}

\begin{remark}\label{rem_daypaper}
In the special case that $S$ is a right cancellative semigroup with a right identity element, our operator $M$ is essentially the same as
the convolution operator $\circ \varphi$ studied by Day \cite{Day64} where $\varphi \in \ell_1(S)$ is the natural
probability mass function induced by the choice of generators:
$$\varphi(s) \ = \ \frac{| \lbrace x \in X \mid \ol{x} = s \rbrace |}{|X|}.$$
Indeed, for any $v \in \ell_2(S)$ and any $s \in S$ we have
\begin{align*}
(vM)(s) \ &= \ \sum_{t \in S} v(t) M_{ts} \ = \ \sum_{t \in S} v(t) \frac{ | \lbrace x \in X \mid t \ol{x} = s \rbrace |}{|X|} \\
&= \ \sum_{t \in S} v(t) \sum_{x \in X, t\ol{x} = s} \frac{1}{|X|} \ = \ \sum_{t \in S} v(t) \sum_{x \in S, tx = s} \varphi(x) \\
&= \ \sum_{t,x \in S, tx = s} v(t) \varphi(x) \ = \ (v \circ \varphi)(s).
\end{align*}
\end{remark}

\begin{remark}
The converse of the final part of Proposition~\ref{prop_operatorcondition} does not hold in general: an operator corresponding to a non-right-cancellative monoid could have norm $1$ or less. This is because the ``concentration'' effect resulting from multiplying by a non-right-cancellative generator can be more than compensated for by ``diffusion'' caused by other generators, as the following example shows.
\end{remark}

\begin{example}
Let $S = \langle a, b, c \mid ac = bc \rangle$, with the obvious choice of generators. Clearly $S$ is not right cancellative. Let $M$ be the corresponding right random walk transition matrix. It is easy to see that any element $s \in S$ can be
written in exactly one of the forms (1) $s = p_s a$ where $p_s \in S$, (2) $s = q_s b$ where $q_s \in S$, (3) $s = r_sac = r_sbc$ where $r_s \in S$
or (4) $s = t_s cc$ where $t_s \in S$. For $i \in \lbrace 1,2,3,4 \rbrace$ we write $S_i$ for the set of elements of $S$ which can be written in form (i), so that $S$ is the disjoint union of $S_1$, $S_2$, $S_3$ and $S_4$. Notice that the choice of $p_s$, $q_s$, $r_s$ or $t_s$ as appropriate is unique.

Now let $v \in \ell_2(S)$. It follows easily from the definition of the operator $M$ that 
$$(vM)(s) \ = \ \begin{cases}
\frac{1}{3} v(p_s) & \textrm{ if } s \in S_1; \\
\frac{1}{3} v(q_s) & \textrm{ if } s \in S_2; \\
\frac{1}{3} v(r_s a) + \frac{1}{3} v(r_s b)  & \textrm{ if } s \in S_3; \\
\frac{1}{3} v(t_s c)  & \textrm{ if } s \in S_4.
\end{cases}$$
Thus we have 
\begin{align*}
(|(vM)|_2)^2 \ = \ &\left( \sum_{s \in S_1} \left(\frac{1}{3} v(p_s) \right)^2 \right)  + \left( \sum_{s \in S_2} \left(\frac{1}{3} v(q_s) \right)^2 \right) \\
&+ \left( \sum_{s \in S_3} \left(\frac{1}{3} v(r_s a) + \frac{1}{3} v(r_s b) \right)^2 \right)
+ \left( \sum_{s \in S_4} \left(\frac{1}{3} v(t_s c) \right)^2 \right).
\end{align*}
Since for each element $z$ there is exactly one element $za$ [respectively, $zb$, $zc$], each element occurs as $p_s$ [respectively, $q_s$, $t_s c$] for at most one
$s$. It follows that the first, second and fourth sums are each bounded above by $\frac{1}{9} (|v|_2)^2$. For the remaining (third) sum we have
$$\sum_{s \in S_3} \left( \frac{1}{3} v(r_s a) + \frac{1}{3} v(r_s b) \right)^2 \ \leq \ \sum_{s \in S_3} \frac{1}{9} 2 \left( v(r_s a)^2 + v(r_s b)^2 \right)$$
by application of Lemma~\ref{lemma_squaresum} to each term. Since for each element $z$ there is only one element $zc$, each element occurs as $r_s a$ or $r_s b$ for at most one $s$. Moreover, the same element cannot occur as both $r_s a$ and $r_s b$ (since no relation allows us to change the last letter of a word), and so each element occurs at most once in the above sum. Thus, this sum is less that $\frac{2}{9} (|v|_2)^2$ and hence $(|vM|_2)^2 \leq \frac{5}{9} (|v|_2)^2$. Since $v$ was a general element of $\ell_2(S)$ the norm (and hence also the spectral radius) of $M$ cannot exceed $\sqrt{\frac{5}{9}} = \frac{\sqrt{5}}{3} \approx 0.745$.
\end{example}

\begin{remark}
A semigroup $S$ is said to have \textit{infinite right indegree} if there are elements $x$ and $s$ such that $tx = s$ for infinitely many
$t$. This will happen, for example, if $S$ is infinite and has a right zero element. In this case, $vM$ is not always defined for $v \in \ell_2(S^1)$.
Indeed, if $S$ has infinite indegree then it is easy to see that we can choose $s$ and $x$ as in the definition, with $x \in X$ a generator. Now
choose $v \in \ell_2(S^1) \setminus \ell_1(S^1)$ supported on the elements $t$ satisfying $tx = s$, and it is clear that $(vM)(s)$ is undefined.
This contrasts with the case of \textit{finite but unbounded} indegree, where for all $v \in \ell_2(S^1)$ the product $vM$ will be defined, but
may not be itself in $\ell_2(S^1)$.
\end{remark}

Our next result says that, just as for groups, maximal cogrowth tells us something about the spectral radius of the corresponding random walk operator.

\begin{theorem}\label{thm_kestenspectralradius}
If a semigroup $S$ has maximal global cogrowth with respect to some finite choice of generators
then the associated right (or left) random walk Markov operator on $\ell_2(S^1)$ is either not defined or has spectral radius at least $1$.
\end{theorem}
\begin{proof}
Suppose the right random walk operator $M$ is defined on $\ell_2(S^1)$.
If we let $\chi_1 \in \ell_2(S^1)$ be once again the characteristic function of the singleton set $\lbrace 1 \rbrace \subseteq S^1$
then for each $n$, combining the definition of the operator norm with Proposition~\ref{prop_cogrowthrandomwalk} yields
$$||M^n||_2 \ \geq \ \frac{|\chi_1 M|_2}{|\chi_1|_2} \ = \ |\chi_1 M|_2 \ = \ \sqrt{\frac{\gamma'(2n)}{|X|^{2n}}}$$
Since $S$ has maximal cogrowth, by Lemma~\ref{lemma_gammaprime} we have
$$\limsup_{n \to \infty} \left( \gamma'(n) \right)^{1/n} \ = \ |X|$$
which means that for any $\beta < 1$ there are infinitely many $k$ such that $\gamma'(k) \geq (\beta |X|)^k$. Since
$\gamma'(k) = 0$ for $k$ odd, this means there are infinitely many $n$ for which $\gamma'(2n) \geq (\beta |X|)^{2n}$,
and hence for which
$$||M^n||_2 \ \geq \ \sqrt{\frac{\gamma'(2n)}{|X|^{2n}}} \ \geq \ \sqrt{\frac{(\beta|X|)^{2n}}{|X|^{2n}}} \ = \ \beta^n.$$
Thus, $M$ has spectral radius at least $\beta$ and since $\beta$ may be chosen arbitrarily close to $1$,
$M$ has spectral radius at least $1$.
\end{proof}

Recall that a semigroup $S$ satisfies \textit{the right F\o lner condition} (\textit{right FC})
if for every finite subset $H$ of $S$ and every $\epsilon > 0$, there
is a finite non-empty subset $F$ of $S$ with $|Fs \setminus F| \leq \epsilon |F|$ for
all $s \in H$.  It satisfies the \textit{right strong F\o lner condition} (\textit{right SFC})
if for every finite subset $H$ of $S$ and every $\epsilon > 0$, there
is a finite non-empty subset $F$ of $S$ with $|F \setminus Fs| \leq \epsilon |F|$ for
all $s \in H$.  There is an obvious dual \textit{left F\o lner condition} (\textit{left FC}) and
\textit{left strong F\o lner condition} (\textit{left SFC}). In groups it is well
known that all of these conditions are equivalent to amenability; in semigroups right FC and
right SFC are respectively necessary and sufficient conditions for right amenability \cite{Paterson88}. It transpires
that right FC is sufficient for the right random walk Markov operator (with respect to
every finite choice of generators) to have spectral radius at least $1$:

\begin{theorem}\label{thm_fcimpspecrad}
If a semigroup $S$ satisfies the right [left] F\o lner Condition then, for any finite choice of generators, the associated right [left] random walk Markov operator on $\ell_2(S^1)$ is either not defined or has spectral radius at least $1$.
\end{theorem}

\begin{proof}
Suppose $S$ satisfies the right F\o lner Condition and that the right random walk 
 operator $M$ is defined. Let $n \in \mathbb{N}$ and $\epsilon > 0$, and consider the operator $M^n$. Let $K$
be the (finite) set of all elements in $S$ represented by words of length $n$ or less. By the F\o lner condition, we may choose a set $F$ such that for each $k \in K$,
$|Fk \setminus F| \leq \frac{\epsilon}{|K|} |F|$. Since $FK \subseteq F \cup \bigcup_{k \in K} Fk \setminus F$, it follows that
$$|FK| \ \leq \ |F| + |K| \frac{\epsilon}{|K|} |F| \ = \ (1+\epsilon) |F|.$$
Let $v \in \ell_2(S^1)$
be the uniform probability mass function on the finite set $F$. A simple calculation shows that $|v|_2 = |F|^{-\frac{1}{2}}$.

Consider now $v M^n$. By the definition of $K$, the only elements of $S$ which can be reached from elements of $F$ by right-multiplication by a sequence of $n$ or
fewer generators are the elements of $FK$. It follows that $v M^n$ is supported only on the set $FK$, which has cardinality at most $(1+\epsilon) |F|$. The $2$-norm
of $v M^n$ therefore cannot be less than the $2$-norm of the uniform probability mass function on $FK$, which is $|FK|^{-\frac{1}{2}}$, so we have
$$|v M^n|_2 \ \geq \ \sqrt{\frac{1}{|FK|}} \ \geq \ \sqrt{\frac{1}{(1+\epsilon) |F|}}$$
and hence
$$||M^n||_2 \ \geq \ \frac{|v M^n|_2}{|v|^2} \ \geq \ \sqrt{\frac{|F|}{(1+\epsilon)|F|}} \ = \ \frac{1}{\sqrt{1+\epsilon}}.$$
Since can be $\epsilon$ chosen arbitrarily small, it follows that $||M^n||_2 \geq 1$ for all $n$, so that $M$ has
spectral radius at least $1$.
\end{proof}

\section{Cogrowth and amenability}\label{sec_amenable}

In this section we consider the relationship between amenability and cogrowth for finitely generated semigroups. We
recall some relevant definitions: for a more detailed introduction the reader is directed to \cite{Paterson88}.

A semigroup $S$ is called \textit{left amenable} if there is a mean
on $l_\infty(S)$ which is invariant under the natural left action
of $S$ on the dual space $l_\infty(S)'$ \cite[Section 0.18]{Paterson88}.
Equivalently \cite[Problem 0.32]{Paterson88}, $S$ is left amenable if
it admits a finitely additive probability measure $\mu$,
defined on all the subsets of $S$, which is \textit{left invariant},
in the sense
that $\mu(a^{-1} X) = \mu(X)$ for all $X \subseteq S$ and $a \in S$. Here,
$a^{-1} X$ denotes the set $\lbrace s \in S \mid as \in X \rbrace$. (Note that
left invariance is strictly weaker than requiring
$\mu(aX) = \mu(X)$ for all $X \subseteq S$ and $a \in S$.)
We also mention an important structural property of semigroups: a semigroup is called \textit{left reversible}
if it does not admit disjoint right ideals; left reversibility is a necessary precondition for left amenability \cite[Proposition~1.23]{Paterson88}.

There are obvious dual notions of \textit{right} invariance, \textit{right} amenability and \textit{right} reversibility. For inverse semigroups (in particular for groups) left and right amenability coincide; a left/right amenable inverse semigroup or group is simply called \textit{amenable} and in fact always admits a measure which is simultaneously left and right invariant. Where a distinction is necessary, we work chiefly 
with left amenability, for consistency with the standard text \cite{Paterson88} and most of the subsequent literature. Of course,
all results admit left/right duals; in some cases we make these explicit, in particular where this gives a  clearer relationship with Section~\ref{sec_randomwalk} above. 

Kesten's Theorem (\cite{Kesten59b}, also stated as Theorem~\ref{thm_kesten} above) can be interpreted as saying that for finitely generated \textit{groups}, amenability is equivalent
to the associated (right or left) random walk Markov operator having spectral radius $1$. In a group, amenability is equivalent to right amenability, the operator is always defined and $1$ is the maximum possible value its spectral radius (by Proposition~\ref{prop_operatorcondition} for example). The theorem can therefore (somewhat vacuously but helpfully in our context --- see Remark~\ref{remark_normintuition}) be rephrased as saying that
right amenability is equivalent to the right random walk operator being undefined or having spectral radius \textit{at least} $1$. Our results so far,
combined with known results about amenability, imply that one implication of this statement holds for semigroups in complete generality:
\begin{theorem}\label{thm_amenableimpspecrad}
If $S$ is a left [right] amenable semigroup then the Markov operator on $\ell_2(S^1)$ of a left [right] random walk on $S$ is either not defined or has spectral radius at least $1$.
\end{theorem}
\begin{proof}
If $S$ is left [right] amenable then by \cite[Proposition 4.9]{Paterson88} or its dual, $S$ satisfies the left [right] F\o lner condition, so the claim follows from Theorem~\ref{thm_fcimpspecrad} and its dual.
\end{proof}

\begin{remark}\label{rem_amenablespecrad}
We note that the converse to Theorem~\ref{thm_amenableimpspecrad} cannot hold for semigroups
in general; indeed if it did we would be able to deduce by combining with Theorem~\ref{thm_fcimpspecrad} that the right F\o lner condition suffices for right amenability, and this is known not to be
the case (see for example \cite[Section~4.22]{Paterson88}).
It remains open whether a converse may hold in well-behaved classes of semigroups, such as left or right cancellative semigroups. Probably the closest thing in the literature is the implication $(e) \implies (a)$ of \cite[Theorem~1]{Day64}: applied to our situation through the translation described in Remark~\ref{rem_daypaper}, this implies that if $S$ is right cancellative with a right identity element and the random walk Markov operator corresponding to \textbf{every} finite choice of generators has norm $1$, then $S$ is right amenable.
\end{remark}

It is a well-known theorem, due separately to Grigorchuk \cite{Grigorchuk80} and Cohen \cite{Cohen82},
that for finitely generated \textit{groups}
considered with symmetric generating sets, amenability is exactly characterised by the maximality
of the \textit{group} cogrowth rate with respect to a symmetric generating set. (Because
only reduced words are considered in the group cogrowth, ``maximal'' in this context means $1$ less than the
size of the generating set). This result can also be rephrased in terms of (either local or
global) semigroup cogrowth for groups:

\begin{theorem}\label{thm_groupamen}
Let $G$ be a group with a finite symmetric choice of generators $X$. Then the following are equivalent:
\begin{itemize}
\item[(i)] $G$ is amenable;
\item[(ii)] $\gamma^{G,X} = |X|$;
\item[(iii)] $\lambda^{G,X}_g = |X|$ for some $g \in G$;
\item[(iv)] $\lambda^{G,X}_g = |X|$ for all $g \in G$.
\end{itemize}
\end{theorem}
\begin{proof}
The equivalence of (ii), (iii) and (iv) follows from Proposition~\ref{prop_localglobalgroup}.
By \cite[Theorem~1]{Elder12} (which is a consequence of the Grigorchuk-Cohen cogrowth theorem
together with a result of Kouksov \cite{Kouksov98}) the group $G$ is amenable if and only if, in our
language, $\lambda^{G,X}_1 = |X|$; since the latter condition implies (iii) and is implied by
(iv), this completes the proof.
\end{proof}

\begin{remark}
The implication (ii) implies (i) of Theorem~\ref{thm_groupamen} can also be deduced more directly from Kesten's Theorem and our results above, without passing through the intervening machinery of the Grigorchuk-Cohen cogrowth theorem. Indeed, if $G$ is group with maximal global cogrowth then the operator $M$ is defined on $\ell_2(G)$ (by Proposition~\ref{prop_operatorcondition}) and has spectral
radius at most $1$ (by Proposition~\ref{prop_operatorcondition}) and at least $1$ (by Theorem~\ref{thm_kestenspectralradius}); hence by Kesten's Theoerm (Theorem~\ref{thm_kesten}) $G$ is amenable. To establish the converse implication in a similar direct way, we would need to prove directly that spectral radius $1$ implies maximum global cogrowth for groups, which does not seem
to be so straightforward.
\end{remark}

When groups are viewed as examples of semigroups, the requirement in Theorem~\ref{thm_groupamen} that the choice of generators be
symmetric is rather unsatisfactory, so it is natural ask whether the symmetry requirement can be dropped. The following results will allow us to
partly answer this question.

\begin{theorem}\label{thm_cogrowthamenable}
Let $S$ be a left reversible, finitely generated semigroup with
maximal global cogrowth, and such that its maximum right cancellative quotient
has a minimal ideal. Then $S$ is left amenable.
\end{theorem}
\begin{proof}
Let $T$ be the maximum right cancellative quotient of $S$.
By Proposition~\ref{prop_globalquotient}, $T$ also has maximal
global cogrowth with respect to the corresponding choice of
generators. Moreover, it is easy to see that $T$ is left reversible.

Let $I$ be the minimal ideal of $T$. We claim that $I$ is a group.
Indeed, $I$ is simple and right cancellative, so it must be a completely
simple semigroup with a single $\GreenL$-class. Moreover, it is easily
seen that $I$ is left reversible, which means it has only one $\GreenR$-class.
Thus, $I$ is a group.

In particular, $I$ is an ideal which is a monoid, and so by Lemma~\ref{lemma_monoidideal} $I$ has maximal global cogrowth
with respect to some choice of generators, and by Theorem~\ref{thm_addidentity} we may assume this choice of generators
contains a representative for the identity element of $I$.

Also, $I$ is left cancellative, so we may apply the dual to Theorem~\ref{thm_kestenspectralradius}, which tells us
that the corresponding left random walk operator $M$
on $\ell_2(I)$
has spectral radius $1$, and hence operator norm $1$. Thus, through the dual to the
translation of notation described in Remark~\ref{rem_daypaper}, the situation
satisfies condition (e'') of the dual to \cite[Corollary to Theorem 4]{Day64}, from
which it follows
that $I$ is an amenable group.
Now by \cite[remarks following Corollary 1.22]{Paterson88} $T$ is left amenable, so by
\cite[Proposition~1.25]{Paterson88} $S$ itself is left amenable.
\end{proof}

\begin{corollary}\label{cor_inverse}
Every finitely generated group or inverse semigroup with maximal global cogrowth is amenable.
\end{corollary}
\begin{proof}
Every inverse semigroup (and hence also every group) is left reversible and has maximum right cancellative quotient which is a group.
Since every group has a minimal ideal, Theorem~\ref{thm_cogrowthamenable} therefore tells us that an inverse semigroup
or group with maximal global cogrowth must be left amenable, which for these semigroups is the same as being amenable.
\end{proof}
We note that Corollary~\ref{cor_inverse} in the group case improves upon the corresponding implication of Kesten's Theorem, because it doesn't
require the choice of generators to be symmetric.

\begin{theorem}\label{thm_localimpamenable}
Every left reversible finitely generated semigroup of maximal local cogrowth is left amenable.
\end{theorem}
\begin{proof}
Suppose $S$ has maximal local cogrowth at an element $s \in S$. Then by Proposition~\ref{prop_localvsglobal} it has maximum global cogrowth. Now since $S$ is left reversible it has a maximum right cancellative quotient; call this $T$. Let $t$ be the image in $T$ of the element $s$. By Proposition~\ref{prop_localquotient}, $T$ has a choice of generators with respect to which $t$ maximal local cogrowth. Now by Proposition~\ref{prop_localminimalideal}, $T$ has a minimal ideal, so the claim follows from Theorem~\ref{thm_cogrowthamenable}.
\end{proof}

\begin{remark}\label{remark_gerl}
Theorem~\ref{thm_localimpamenable} is related to a theorem of Gerl \cite[Theorem~3]{Gerl73}.
He considers probability measures $P$ on a semigroup $S$ such that the support of $P$ generates $S$,
defines a convolution product on measures by
$$(P*Q)(s) = \sum_{s = s_1 s_2} P(s_1) Q(s_2)$$
and considers the powers $P^{(n)}$ of $P$
under this operation. His theorem states that if $S$ is left cancellative with a left
unit element and 
$$\limsup_{n \to \infty} \left( P^{(n)}(s) \right)^{1/n} \ = \ 1$$
then $S$ is left amenable.

Returning to our usual setting, with $X$ a finite choice of generators for $S$, we may define a (finitely supported) measure $P$ on $S$ where each
element is weighted in proportion to the number of generators in $X$ representing it. Then it is easy to see that
$$P^{(n)}(s) = \frac{\lambda_s(n)}{|X|^n}$$
for every $n$. In particular, if $S$ has maximal local cogrowth at $s$ (and satisfies the structural conditions of being left cancellative with a left unit) then the hypotheses of Gerl's theorem are satisfied and so we obtain an alternative way to deduce that $S$ is left amenable.

Gerl's theorem is more general than Theorem~\ref{thm_localimpamenable} in that it treats a general probability measure (possibly with infinite support) and hence has potential for application to non-finitely generated semigroups, but more specific in that it applies only when $S$ is left cancellative with a left unit. Gerl gives counterexamples \cite[Remark 3]{Gerl73} to show that the latter structural assumptions cannot be completely removed, but the obstruction to left amenability in these cases is always the absence of left reversibility. We conjecture
that the left cancellativity and left unit hypotheses in Gerl's theorem can be replaced with left reversibility:
\end{remark}

\begin{conjecture}
Let $P$ be a probability measure on a (not necessarily finitely generated) left reversible semigroup $S$, such that the support of $P$ generates $S$.  Define the convolution power $P^{(n)}$ as in Remark~\ref{remark_gerl}. If 
$$\limsup_{n \to \infty} \left( P^{(n)}(s) \right)^{1/n} \ = \ 1$$
then $S$ is left amenable.
\end{conjecture}

We now return to the case of a finitely generated group $G$, and the question of whether the characterisations of amenability given by Theorem~\ref{thm_groupamen} remain equivalent if the symmetry assumption on the choice of generators is dropped. If $G$ has maximal global cogrowth then by Corollary~\ref{cor_inverse} it must be amenable. Similarly, if $G$ maximal local
cogrowth then by Theorem~\ref{thm_localimpamenable} (or by Proposition~\ref{prop_localvsglobal} and Corollary~\ref{cor_inverse}) it is amenable.

The converse fails for \textit{local} cogrowth:  indeed, in Example~\ref{example_groupnonsymmetric_local} above we saw a $3$-generated amenable group ($\mathbb{Z}$ with an asymmetric choice of generators) with local cogrowth $2 \sqrt{2} < 3$ at every element.  However, by Corollary~\ref{cor_commutativeglobal}, the \textit{global} cogrowth of $\mathbb{Z}$ is maximal with respect to any (symmetric or asymmetric) choice of generators, so this cannot be a counterexample to the converse for global cogrowth (and nor can any other group of subexponential growth).
\begin{question}\label{q_groupamenableimpglobal}
Does there exist an amenable group $G$ (necessarily of exponential growth) with a choice of generators $X$ (necessarily asymmetric)
such that $\gamma^{G,X} < |X|$?
\end{question}
While we do not know whether or not there is an amenable group without maximal global cogrowth, we do expect
that there are semigroups and monoids with the corresponding property:
\begin{conjecture}\label{conjecture_amenablecogrowth}
There is a left reversible finitely generated monoid which is left amenable but does not have maximal global cogrowth.
\end{conjecture}

Recall from \cite[Section 2.4]{GrayKambitesAmenable} that a semigroup $S$ satisfies the \emph{Klawe condition} if whenever $s$, $x$ and $y$ in $S$ are such that $sx = sy$, there exists $t \in S$ so that $xt = yt$. This condition first arose implicitly, without being given a name, in the work of Klawe \cite{Klawe77}. It is immediate from the definition that every left cancellative semigroup satisfies the Klawe condition. It follows from \cite[Propositions~2.2 and 2.3]{GrayKambitesAmenable} that every left reversible semigroup in which every ideal contains an idempotent also satisfies the Klawe condition; this includes every group, inverse semigroup, left reversible regular semigroup, left reversible finite semigroup, left reversible compact left or right topological semigroup, and semigroup with a left, right or two-sided zero element.

\begin{theorem}\label{thm_cogrowthamenable_Klawe}
Let $S$ be a left reversible, finitely generated semigroup such that $S^1$ satisfies the Klawe condition.
If $S$ has maximal global cogrowth then $S$ is left amenable. 
\end{theorem}
\begin{proof}
By Theorem~\ref{thm_supersets} the monoid $S^1$ has maximal cogrowth with respect to some choice of generators; call it $A$.
Let $T$ be the maximum right cancellative quotient of $S^1$. By Proposition~\ref{prop_globalquotient}, $T$ has maximal global cogrowth with respect to $A$, and by \cite[Proposition~2.4]{GrayKambitesAmenable}, $T$ is a left cancellative monoid. By Theorem~\ref{thm_supersets}, for any finite subset $K$ of $T$ there exists a finite generating set $B_K$ for $T$ such that, $B_K$ contains the identity element of the monoid $T$, $K$ is a subset of $B_K$, and $T$ has maximal cogrowth with respect to $B_K$. By Theorem~\ref{thm_kestenspectralradius} it follows that (with respect to this generating set $B_K$) the corresponding left Markov operator $M$ on $\ell_2(T)$ has spectral radius $1$ and hence operator norm $1$.

Thus the dual of condition (e) in \cite[Theorem 1]{Day64} is satisfied. (The general translation into the language of \cite{Day64} is as described in Remark~\ref{rem_daypaper} above; the set $K$ corresponds to $\xi$ from \cite[Theorem~1]{Day64}, the set $P_{\varphi} \cap U$ is not empty since $B$ contains the identity element, and we take $p=2$ since we are working with the $\ell_2$-norm of the operator $M$.)

Applying that result we conclude that $T$ is left amenable, and hence by \cite[Proposition~1.25]{Paterson88} $S^1$ is left amenable. Since $S$ is a left ideal in $S^1$, it follows by \cite[Corollary 1.22]{Paterson88} that $S$ is left amenable.
\end{proof}

\begin{remark}\label{remark_s1klawe}
The hypothesis in Theorem~\ref{thm_cogrowthamenable_Klawe} that $S^1$ satisfies the Klawe condition is slightly
unsatisfactory, but in many cases of interest it coincides with the (more natural) condition that $S$ satisfies the Klawe condition. Indeed, if $S$ satisfies the Klawe condition, the only possible obstruction to $S^1$ satisfying the Klawe condition is the existence of elements $s,x \in S$ such that $sx = s$ but there is no $t \in S$ with $xt = t$.
\end{remark}

\begin{corollary}\label{cor_regular_left_reversible}
Every left reversible finitely generated regular semigroup with maximal global cogrowth is left amenable.
\end{corollary}
\begin{proof}
If $S$ is left reversible and regular then clearly so is $S^1$, so by \cite[Propositions~2.2 and 2.3]{GrayKambitesAmenable}
the monoid $S^1$ satisfies the Klawe condition and the result follows from Theorem~\ref{thm_cogrowthamenable_Klawe}.
\end{proof}

\begin{corollary}\label{cor_leftcancellative}
Every left reversible left cancellative finitely generated monoid with maximal global cogrowth is left amenable. 
\end{corollary}

\begin{remark}\label{rem_cogrowthamenable}
Theorems~\ref{thm_cogrowthamenable}, \ref{thm_localimpamenable} and \ref{thm_cogrowthamenable_Klawe} and Corollaries~\ref{cor_inverse}, \ref{cor_regular_left_reversible} and \ref{cor_leftcancellative} all give sufficient structural conditions for maximal local or global cogrowth to imply left amenability. In all cases
left reversibility is (explicitly or implicitly) among the conditions, and indeed this condition is necessary, since (for example) all finite semigroups have maximal global cogrowth (by Corollary~\ref{cor_finitekesten}) and maximal local cogrowth (by Proposition~\ref{prop_minimalideal}), but left reversibility is necessary for left amenability \cite[Proposition~1.23]{Paterson88}. In fact, by the latter observation and Theorem~\ref{thm_localimpamenable}, a semigroup of maximal local cogrowth is left amenable if and only if it is left reversible. We do not know if the same holds for maximal global cogrowth:
\end{remark}
\begin{question}\label{qn_cogrowthamenable}
Is a left reversible finitely generated semigroup of maximal global cogrowth necessary left amenable?
\end{question}

\section{Conclusions and open problems}

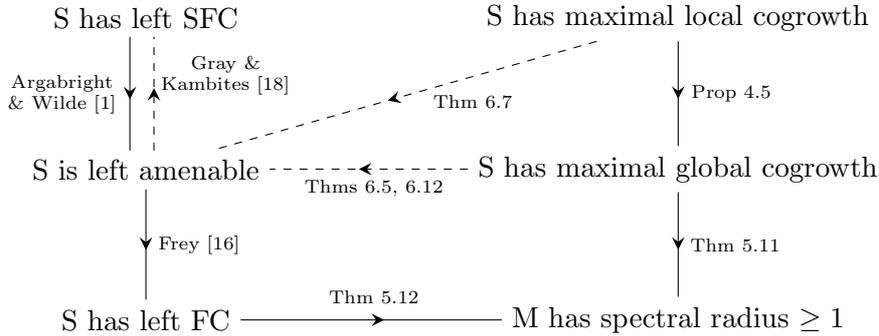
\begin{figure}
\def\x{0.55}
\begin{center}
\begin{tikzpicture}[decoration={ 
markings,
mark=
at position \x
with 
{ 
\arrow[scale=1.5]{stealth} 
} 
} 
]
\tikzstyle{vertex}=[circle,draw=black, fill=white, inner sep = 0.75mm]
\tikzset{
    punkt/.style={
           rectangle,
           rounded corners,
           draw=black, very thick,
           text width=8em,
           minimum height=2em,
           text centered}
    pil/.style={
           ->,
           thick,
           shorten <=2pt,
           shorten >=2pt,}           
}           
%
%
\node (A) at (3,0.3) {\tiny{Thm~\ref{thm_fcimpspecrad}}};
\node (D) at (0.7,1) {\tiny{Frey \cite{Frey60}}};
\node (E) at (7.75,1) {\tiny{Thm~\ref{thm_kestenspectralradius}}};
\node (F) at (7.7,3){\tiny{Prop~\ref{prop_localvsglobal}}};
\node (G) at (4.3,2.9) {\tiny{Thm~\ref{thm_localimpamenable}}};
\node (G) at (3,1.75) {\tiny{Thms~\ref{thm_cogrowthamenable},~\ref{thm_cogrowthamenable_Klawe}}};
\node  (H) at (-1.1,3) {\tiny{\begin{tabular}{c} Argabright \\ \& Wilde \cite{Argabright67} \end{tabular}}};
\node (I) at (1.05,3.25) {\tiny{\begin{tabular}{c} Gray \& \\ Kambites \cite{GrayKambitesAmenable} \end{tabular}}};
%
%
\node[punkt] (rSFC) at (0,4) {S has left SFC};
\node[punkt] (rA) at (0,2) {S is left amenable};
\node[punkt] (rFC) at (0,0) {S has left FC};
\node[punkt] (mLC) at (7,4){S has maximal local cogrowth};
\node[punkt] (mGC) at (7,2) {S has maximal global cogrowth};
\node[punkt] (rSpecR) at (7,0) {M has spectral radius $\geq 1$};
%
%
%
%
%
\arrowdraw{mGC}{rSpecR}
\darrowdraw{mLC}{rA}
\arrowdraw{mLC}{mGC}
\arrowdraw{rFC}{rSpecR}
\arrowdraw{rA}{rFC}
\arrowdraw{[xshift=-6pt,yshift=-7pt]rSFC.center}{[xshift=-6pt,yshift=7pt]rA.center}
%
%
\darrowdraw{mGC}{rA}
\darrowdraw{[xshift=3pt,yshift=7pt]rA.center}{[xshift=3pt,yshift=-7pt]rSFC.center}

\end{tikzpicture}
\end{center}
\caption{Relationships between various properties of a finitely generated semigroup $S$ and of the corresponding left random walk Markov operator $M$.  
} 
\label{fig_relationships}
\end{figure}

Figure~\ref{fig_relationships} summarises some of the main results of this article and also some relevant results of our previous article \cite{GrayKambitesAmenable} and older work of Frey \cite{Frey60} and Argabright and Wilde \cite{Argabright67}. It shows seven different conditions
which for groups with symmetric generating sets are known to coincide, essentially by the work of F\o lner \cite{Folner55}, Kesten \cite{Kesten59b},
Grigorchuk \cite{Grigorchuk80} and Cohen \cite{Cohen82}. Solid arrows represent implications which hold for finitely generated semigroups in complete generality. Dashed arrows represent implications known to hold under certain mild structural assumptions on the semigroup.

All of the dashed arrow implications hold 
for left reversible, finitely generated monoids satisfying the Klawe condition, and hence for all groups (with not necessarily symmetric generating sets), inverse monoids, left reversible regular monoids, left reversible left cancellative and left reversible near left cancellative monoids. They also hold for left reversible regular semigroups. Of the missing implications, we know maximal local cogrowth is not equivalent to left amenability (even for groups with asymmetric choice of generators --- see Example~\ref{example_groupnonsymmetric_local}). Within this class of monoids, we do not know whether left amenability implies maximal global cogrowth (or equivalently, left SFC), even for groups with asymmetric choice of generators
(Question~\ref{q_groupamenableimpglobal}). Since global cogrowth is not ``sided'', a positive answer would imply that for semigroups which are
left and right reversible and satisfy the Klawe condition and its dual, left and right amenability coincide.
We also do not know whether the left Markov operator having spectral radius (or norm) at least $1$ suffices for left amenability in this
setting --- see
Remark~\ref{rem_amenablespecrad}; if so this
would imply that (within the given class of monoids) left FC coincides with left amenability and left SFC.

In the wider class of left reversible finitely generated semigroups (without the Klawe condition), we know that maximal local cogrowth implies left amenability (Theorem~\ref{thm_localimpamenable}), but we do not know if maximal global cogrowth, or a Markov operator with norm and/or spectral radius at least $1$ suffice for left amenability, or even for left FC. Left amenability certainly does not imply maximal local cogrowth (Example~\ref{example_freecommutativelocal}), but we do not know whether it suffices for maximal global cogrowth; we conjecture that it does not (Conjecture~\ref{conjecture_amenablecogrowth}).

For finitely generated semigroups in absolute generality, it is well known that left FC does not imply left amenability (\cite{Frey60}, or see \cite[Section~4.22]{Paterson88} for a more accessible reference) and that left amenability does not imply left SFC (due to an example
of \cite{Takahashi03} which is a refinement to the finitely generated case of a result of Klawe \cite{Klawe77}). We have seen that 
the existence of an element of maximal local cogrowth does not imply left amenability (see Remark~\ref{rem_cogrowthamenable} above), and
therefore nor does maximal global cogrowth, or the associated random walk Markov operator having a norm or spectral radius of at least $1$.

\bibliographystyle{plain}

\def\cprime{$'$} \def\cprime{$'$}

\end{document}